
\input amstex
\input amsppt.sty
\magnification=\magstep1
\hsize = 6.5 truein
\vsize = 9 truein

\NoRunningHeads
\NoBlackBoxes
\TagsAsMath

\def\label#1{\par%
        \hangafter 1%
        \hangindent .75 in%
        \noindent%
        \hbox to .75 in{#1\hfill}%
        \ignorespaces%
        }

\newskip\sectionskipamount
\sectionskipamount = 24pt plus 8pt minus 8pt
\def\sectionskip{\vskip\sectionskipamount}
\define\sectionbreak{%
        \par \ifdim\lastskip<\sectionskipamount
        \removelastskip    \penalty-2000 \sectionskip \fi}
\define\section#1{%
        \sectionbreak 
        \subheading{#1}%
        \bigskip
        }

\redefine\qed{{\unskip\nobreak\hfil\penalty50\hskip2em\vadjust{}\nobreak\hfil
    $\square$\parfillskip=0pt\finalhyphendemerits=0\par}}

\define\op#1{\operatorname{\fam=0\tenrm{#1}}} 
\let \< = \langle
\let \> = \rangle
\define \a {\alpha}
\define \e {\varepsilon}
\define	\x  {\times}
	\redefine	\l		{\lambda}
	\redefine	\L		{\Lambda}

\topmatter
\title Integrable and Proper Actions on $C^*$-Algebras, \\
and Square-Integrable Representations of Groups
\endtitle
\author Marc A. Rieffel
\endauthor
\thanks The research was supported in part by National Science Foundation grant
DMS96-13833.
\endthanks
\abstract We propose a definition of what should be meant by a {\it proper}
action of a locally compact group on a $C^*$-algebra. We show that when
the $C^*$-algebra is commutative this definition exactly captures the 
usual notion of a proper action on a locally compact space. We then
discuss how one might define a {\it generalized fixed-point algebra}. The
goal is to show that the generalized fixed-point algebra
is strongly Morita equivalent to an ideal in the crossed product algebra, as
happens in the commutative case. We show that one candidate 
gives the desired algebra when the $C^*$-algebra is
commutative. But very recently Exel has shown that this candidate
is too big in general. Finally, we consider in detail the application
of these ideas to actions of a locally compact group on the algebra of
compact operators (necessarily coming {}from unitary representations), and
show that this gives an attractive view of the subject of square-integrable
representations.
\endabstract
\address\hskip-\parindent
        Department of Mathematics \newline University of California
\newline Berkeley, California 94720-3840
\endaddress
\date  June 1999
\enddate
\email rieffel\@math.berkeley.edu
\endemail
\endtopmatter

\baselineskip=14 pt
\document

There is a variety of situations in which actions of locally compact groups on
non-commutative $C^*$-algebras appear to be ``proper'' in a way analogous
to proper
actions of groups on spaces.  See for example 
\cite{OP1, OP2, Ks, Rf7, Ma, Qg, Rf8, QR, Ab, E1, GHT}.  We
propose
here a simple definition which seems to capture this idea reasonably well.  We
indicate a variety of examples, but we only explore two basic ones in detail.
Namely, we show that when the $C^*$-algebra is commutative our definition does
capture exactly the usual definition of a proper action on a space.  Then
we show
that when the $C^*$-algebra is the algebra of compact operators on a
Hilbert space
our definition is very closely related to square-integrable 
representations (not necessarily irreducible) of
groups, and gives an attractive view-point on
this
venerable subject.

Our definition of proper actions is closely related to ideas of ``integrable''
actions which occur in various places, especially in the literature concerning
actions on von~Neumann algebras \cite{CT,Pa,S}.  We give here a definition of
``integrable'' actions for $C^*$-algebras which appears to be the right
analogue of
that for von~Neumann algebras.  We see that every proper action is
integrable, but
not conversely.  But it is useful to see that some of the basic properties of
proper actions come just {}from the fact that they are integrable.

I had earlier given a tentative definition of proper actions \cite{Rf7},
which was
adequate to treat some interesting examples.  But that definition assumed the
existence of a dense subalgebra with certain properties, and so was not
intrinsic.
The definition proposed here is intrinsic, and includes my older
definition.  But
the definition given here must still be viewed as tentative, since I have
not yet
been able to relate it strongly to the crossed product $C^*$-algebra for the
action
in the way done in \cite{Rf7} (and there are many other aspects which also
still need to be explored).  

The main goal is to define a suitable generalized fixed-point algebra
(corresponding to the orbit space of a proper action on a space), which will
in a natural way be strongly Morita equivalent to an ideal in the crossed
product algebra, as done in \cite{Rf7}. In an earlier version of this paper
I had proposed a candidate for this generalized fixed-point algebra, and shown
that it works correctly when the algebra acted upon is Abelian (Theorem 6.5
below). But in a very recent preprint \cite{E2} Exel gives a natural example 
showing that in general the candidate which I had proposed is too big. He
also gives a penetrating analysis of the difficulties involved, already for
the case when the group which acts is Abelian. He does this within the context
of making
strong progress on his project of determining which actions on $C^*$-algebras
are dual actions on Fell bundles. But this leaves unresolved 
the question of whether there is
an intrinsic definition of the generalized fixed-point algebra.

In spite of this unsatisfactory situation, it seems to me worthwhile publishing
what I discovered. There is not much overlap between the very recent paper of Exel
and the present paper, and Exel makes use of some of the examples and
results of the present paper. Also, I needed to sort out some of the issues
discussed here
for use in connection with several of my projects concerning quantization.  
(And even the
classical notion of
proper actions on spaces is one of continuing strong interest \cite{BCH, GHT}.)

Actually, the definition of proper actions which we give here was strongly stimulated
by a slightly earlier paper of Exel \cite{E1} (which 
in turn built on \cite{Rf7}).  In fact our
definition almost appears explicitly in \cite{E1}.  The main difference is
that here
we emphasize the order properties of $C^*$-algebras while in \cite{E1} the
emphasis
is on unconditional integrability (as it is in \cite{E2} also).  
I am very grateful to Exel for quite
helpful comments about all these matters.

We will see that our definition of proper actions is closely related to 
earlier notions of integrable elements discussed in
\cite{Ld, OP1, OP2, Qg, QR}. It is also closely related to the
notion
of $C^*$-valued weights on $C^*$-algebras which was introduced recently by
Kustermans \cite{Ku}, for fairly different reasons involving Haar measures for
quantum groups.  (I thank Kustermans for some helpful comments on this matter.)
One can in turn ask what should be meant by proper actions of
quantum groups.  Integrable actions of a fairly wide class of quantum
groups (namely Kac
algebras) acting on von~Neumann algebras are discussed in section $18.19$ of
\cite{S}.  The action of any compact quantum group on a $C^*$-algebra should be
proper, and
indeed in this case one obtains the kind of relations between the fixed point
algebra and the crossed product algebra which one expects \cite{Ng}.  One
can also
ask about proper actions of groupoids on $C^*$-algebras, extending the
notion of
proper actions of groupoids on spaces given in \cite{Re}.

In section~1 we deal with integrable actions, while in section~2 we discuss
$C^*$-valued weights. Section~3 is concerned with the special case of
algebras of continuous functions on a locally compact space with
values in a $C^*$-algebra. In section~4 we combine the earlier material
to define and discuss proper actions. The functoriality properties 
of the situation are discussed in section~5. We also show there that
the commonly used structure of C*-algebras ``proper over an action
on an ordinary space'' \cite{Ks, GHT} falls not only within our present
context of proper actions, but, even more, within the context of
\cite{Rf7}, where strong Morita equivalence of the generalized fixed-point
algebra with an ideal of the crossed product algebra is established.
In section~6 we discuss how one might define the
generalized fixed-point algebra for our present setting. Section~7 is
devoted to actions on the algebra of compact operators, and
their relation to square-integrable representations. Then in section~8 we
continue that
discussion by considering the orthogonality relations.  Substantial parts of
sections~7 and 8 can be viewed as expository, treating well-known material
{}from a
slightly different point of view.

\section{1.  Integrable Actions}

The material discussed here is very close to material on integrable actions
in the
von~Neumann algebra literature.  See for example definition $2.1$ of \cite{CT},
the introduction to \cite{Ld}, 
\cite{Pa}, and $18.20$ of \cite{S}.  Here we stress the $C^*$-algebra
version of
integrable actions, so that we can contrast it with the notion of proper
actions
which we discuss in the next section.  Since every proper action is integrable,
this section also develops those facts about proper actions which hold because 
they are integrable actions.

Let $G$ be a locally compact group, equipped with a choice of left Haar
measure.
Let $\alpha$ be a (strongly continuous) action of $G$ on a $C^*$-algebra
$A$.  (A
simple but useful example to keep in mind during the following discussion
is the
case of $G = {\Bbb R}$ acting on the one-point compactification, ${\tilde {\Bbb
R}}$, of ${\Bbb R}$ by translation, leaving the point at infinity fixed, and so
acting on $A = C({\tilde {\Bbb R}})$.  But in general we do not assume that
$A$ has
an identity element.) Notice that for given $a\in A$ the function
$x \mapsto \alpha_x(a)$ has constant norm, and so can not be integrable over
$G$ (unless $a=0$) when $G$ is not compact. Nevertheless, our aim is to 
give meaning to
  $$ \int_G \alpha_x (a) dx  ,$$
at least for some actions $\alpha$ and some $a\neq 0$.

It is convenient initially to place this matter in a more general context.
Let $X$ be a locally compact space (e.g. $G$) and fix on it a positive
Radon measure (e.g. Haar measure), to which we will not give a particular
symbol. Consider the C*-algebra $B=C_b(X,A)$ of bounded norm-continuous
functions {}from $X$ to $A$. We can surely integrate functions of compact
support. But there may be other functions, even ones of constant norm,
whose integrals we can make sense of indirectly.

For the case of $G$ and $\alpha$, we identify $a\in A$ with the function
$x \mapsto \alpha_x(a)$ in $C_b(G, A)$. This gives an isometric inclusion
of $A$ as a C*-subalgebra of $C_b(G, A)$ (consisting entirely of functions
of constant norm), whose image we will denote by $A_\alpha$. 
So we see that it may be useful in our more general
case of $X$ to consider eventually various subalgebras of $B$. For example, 
our considerations can be applied to the induced C*-algebras studied
for instance in \cite{QR}. Here one has both an action $\alpha$ on $A$
and an action, $\tau$, on a space $M$, and one considers the subalgebra
of $C_b(M,A)$ consisting of the functions $f$ such that $f(\tau_x^{-1}(m))
=\alpha_x(f(m))$. 

For any positive 
$\lambda \in L^1(X)$ define a linear map, $p_\lambda$, {}from $B$
to $A$ by
  $$p_\lambda(f) = \int f(x)\lambda(x) dx.$$
It is easily seen that $p_\lambda$ is positive, in fact completely
positive \cite{KR2}, and of norm $\|\lambda \|_1$. We would like to have the
flexibility of having $\lambda$ range over characteristic
functions of
compact sets, or over continuous functions which approximate them.  It is thus
convenient for us to set, for use throughout this paper,
$$
{\Cal B} = {\Cal B}(X) = \{\lambda \in L^{\infty}(X): \lambda 
\text{ has compact
support and }
0 \le \lambda \le 1\}.
$$
We note that if $\lambda \in {\Cal B}$ then $\lambda \in L^1(X)$, 
so that $p_\lambda$ is
defined.  Also, ${\Cal B}$ is an upward directed set under the usual
ordering of
functions, and if $\lambda_1, \lambda_2 \in {\Cal B}$ 
with $\lambda_1 \le \lambda_2$, 
then $p_{\lambda_1} \le
p_{\lambda_2}$
for the usual ordering of positive maps.  
Thus $\{p_\lambda\}_{\lambda \in {\Cal
B}}$ is
an increasing net of completely positive maps {}from $B$ into $A$.  Let $f \in
B^+$ ( the positive part of $B$).  Then $\{p_\lambda(f)\}_{\Cal B}$ 
is an increasing net of positive
elements of
$A$.  Even if this net is bounded, we can not expect it to converge in $A$.
But
bounded increasing nets of positive elements do converge (for the strong
operator
topology) if they are in a von~Neumann algebra ($5.1.4$ of \cite{KR2}).
Thus if we
view $p_\lambda$ as taking values in the double-dual (or ``enveloping'')
von~Neumann
algebra, $A''$, of $A$ (see 3.7 of \cite{Pe2}), 
we will have such convergence.
Let us examine this situation a bit more generally.

\definition{1.1 Definition} 
Let $N$ be some
von~Neumann
algebra, and let $P = \{p_{\lambda}\}$ be an increasing net of
positive
maps {}from a C*-algebra $B$ to $N$. 
We say that $b \in B^+$ is $P$-{\it bounded} if
the net
$\{p_{\lambda}(b)\}$ is bounded above.  Let ${\Cal M}_P^+$ denote the set of
$P$-bounded elements.  For each $b \in {\Cal M}_P^+$ let $\varphi_P(b)$
denote the
least upper bound of $\{p_{\lambda}(b)\}$ in $N$.  We call the mapping
$\varphi_P$
{}from ${\Cal M}_P^+$ to $N$ the {\it least upper bound} of the net $P$.
\enddefinition

It is evident that ${\Cal M}_P^+$ is a hereditary cone in $B^+$, and that
$\varphi_P$ is ``linear'' and positive.
Now for any hereditary cone ${\Cal M}^+$ in a $C^*$-algebra $B$ we have the
following
structure.  (See $7.5.2$ of \cite{KR2} or $5.1.2$ of \cite{Pe2}.)  Let
${\Cal N} =
\{b \in B: b^*b \in {\Cal M}^+\}$.  Then ${\Cal N}$ is a left ideal in $B$ (not
necessarily closed).  Let ${\Cal M}$ be the linear span of ${\Cal M}^+$.  Then
${\Cal M} = {\Cal N}^*{\Cal N}$ (linear span of products), 
and ${\Cal M} \cap B^+ = {\Cal M}^+$.  In
particular, ${\Cal M}$ is a hereditary $*$-subalgebra of $B$.  If $\varphi$
is an additive map {}from ${\Cal M}^+$ to $M^+$ for even some $C^*$-algebra $M$, 
and if
$\varphi(ra) = r\varphi(a)$ for $r \in {\Bbb R}^+$ and $a \in {\Cal M}^+$,
then the
usual proof for scalar-valued weights shows that $\varphi$ has a positive
linear
extension to ${\Cal M}$, which is unique. An important slippery point
in connection with all this is that if $a \in {\Cal M}$, it does not
follow in general that $|a| \in {\Cal M}$ (even if $a = a^*$). This
difficulty already occurs with ordinary weights. See the example
following theorem 2.4 of \cite{Pe1}. This makes it awkward to define
an ``$L^1$-norm'' on ${\Cal M}$ using $\varphi$. 
In Theorem 8.9 we will give a
simple explicit example in which this difficulty occurs exactly in
our context, namely for an integrable action on the algebra of
compact operators.

Returning to our situation of an increasing net $P = \{p_{\lambda}\}$ with
least
upper bound $\varphi_P$, we see that $\varphi_P$ extends to the linear
span, ${\Cal
M}_P$, of ${\Cal M}_P^+$.  Then it is clear that for $b \in {\Cal M}_P$ the net
$\{p_{\lambda}(b)\}$ converges strongly to $\varphi_P(b)$.  It is not
difficult to verify that if each
$p_{\lambda}$ is completely positive, then $\varphi_P$ is also.
This
suggests the following definition, where we momentarily allow values in a
$C^*$-algebra rather than a von~Neumann algebra.  This definition is
essentially
$1.1$ of \cite{Ku} with $C = M(A)$.

\definition{1.2 Definition} By a $C^*$-{\it valued weight} on a
$C^*$-algebra $B$
we mean a function, $\varphi$, whose domain is a hereditary cone ${\Cal
M}^+$ in
$B^+$, and whose range is contained in $C^+$ for some 
$C^*$-algebra $C$, such that

\roster
\item"1)" $\varphi(ra) = r\varphi(a)$ for $a \in {\Cal M}^+$ and $r \in
{\Bbb R}^+$,

\item"2)" $\varphi(a+b) = \varphi(a) + \varphi(b)$ for $a,b \in {\Cal M}^+$.
\endroster
We will say that $\varphi$ is {\it completely} positive if
the unique positive extension of $\varphi$ to the linear span,
${\Cal M}$,
of
${\Cal M}^+$ is completely positive, in the sense that for all $n$, if
$(b_{ij})$
is an $n \times n$ matrix of elements of ${\Cal M}$ which is positive as an
element
of $M_n(B)$, then the matrix $(\varphi(b_{ij}))$ is positive as a matrix in
$M_n(C)$.  If the values are in a von~Neumann algebra, we will refer to
$\varphi$
as an {\it operator-valued weight on $B$}.
\enddefinition

In the case in which $\varphi$ comes {}from an increasing net of
positive
maps as above, with values in a von~Neumann algebra $N$, it is natural in
view of
standard definitions in the literature (see $5.1.1$ of \cite{Pe2}), to make the
following definition:

\definition{1.3 Definition} An operator-valued weight $\varphi$ on a
$C^*$-algebra
$B$, with domain ${\Cal M}^+$ and range in $N$, is said to be {\it normal}
if there
is an increasing net $\{p_{\lambda}\}$ of bounded positive
linear maps
{}from $B$ into ${\Cal N}$ such that

\roster
\item"1)" ${\Cal M}^+ = \{b \in B^+: \{p_{\lambda}(a)\} \text{ is bounded
above}\}$,

\item"2)" $\varphi(b) = \text{ l.u.b.}\{p_{\lambda}(b)\}$ for $b \in {\Cal
M}^+$.
\endroster
\enddefinition

We return to the situation in which $B = C_b(X,A)$.

\definition{1.4 Definition} Let $B = C_b(X,A)$, and let $P = \{p_\lambda\}$
be as defined earlier.
The elements in the linear span of the $P$-bounded
elements will be called the {\it order-integrable} elements of $B$.
\enddefinition

This definition is closely related to Exel's definition of pseudo-integrable
elements \cite{E1, E2}, the difference being that we emphasize the order
structure rather than the unconditional integrability.  
It is different {}from the definition given in $7.8.4$ of
\cite{Pe2}. Rather we will see that the latter is very close to our
definition of
{\it proper} actions given in the next section.

By considering the continuity of the $p_{\lambda}$'s we obtain the following
alternative characterization of positive order-integrable elements in
this case:

\proclaim{1.5 Proposition} An element $f \in B^+$ is order-integrable
iff there
is a constant, $k_f$, such that
$$
\|p_\lambda(f)\| \le k_f\|\lambda\|_{\infty}
$$
for every $\lambda \in L^{\infty}(G) \cap L^1(G)$ with $\lambda \ge 0$.  
(Equivalently,
we can
omit the condition $\lambda \ge 0$.)
\endproclaim

With a possible change in the constant $k_f$, we obtain:

\proclaim{1.6 Corollary} For every order-integrable element $f \in B$
there is a
constant, $k_f$, such that
$$
\|p_\lambda(f)\| \le k_f\|\lambda\|_{\infty}
$$
for all $\lambda \in L^1(G) \cap L^{\infty}(G)$.
\endproclaim

\demo{\bf 1.7 Notation} We will denote the hereditary $*$-subalgebra of
order-integrable elements by ${\Cal M}_X$, and the left ideal
$\{a \in A:
a^*a \in {\Cal M}_X\}$ by ${\Cal N}_X$.  We denote the associated
operator-valued weight with values in $A''$, and its unique extension to ${\Cal
M}_X$, by $\varphi_X$.  It is natural to also denote
$\varphi_{X}(f)$ for $f \in {\Cal M}_X$ by
$$
\varphi_X(f) = \int f(x)dx,
$$
as long as the integral is interpreted as simply meaning $\varphi_X(f)$.
\enddemo

We will later find the following fact useful.

\proclaim{1.8 Proposition} Let $f \in {\Cal M}_X$, and let $\omega
\in A'$,
the dual space of $A$.  Then the function $x \mapsto \omega(f(x))$ is
integrable (in the ordinary sense) on $X$, and
$$
\int \omega (f(x))dx = \omega(\varphi_X(f)),
$$
where $\omega$ is viewed as being in the predual of $A''$.
\endproclaim

\demo{Proof} By the definition of ${\Cal M}_X$ and by the standard
decomposition \cite{Pe2} of elements of $A'$ in terms of positive elements, it
suffices to treat the case of positive $a$ and positive $\omega$.  The
function in
question is continuous, so measurable.  For $\lambda \in {\Cal B}$ we have
$$
\int \lambda(x)\omega(f(x))dx = \omega(p_\lambda(f)) \le \|\omega\|k_f,
$$
where $k_f$ is as in Corollary $1.6$.  Then a short argument 
using the monotone convergence
theorem shows that the function is integrable.  The weak-$*$
topology on $A''$ coincides with the ultra-weak operator topology
($3.5.5$-$6$ of
\cite{Pe2}), so the equality must hold. \qed
\enddemo

As indicated above, we are interested in subalgebras of B. The main
definition of this section is:

\definition{1.9 Definition} Let $B=C_b(X,A)$ as above, and let  $C$ be
a C*-subalgebra of $B$. We will say that $C$ is an {\it integrable} 
subalgebra if $C \cap {\Cal M}_X$ is dense in $C$.
\enddefinition

We remark that this definition can even be applied to operator systems,
i.e. self-adjoint subspaces, and might eventually be useful there.
Our main application of this definition is to the case of an 
action $\alpha$ of $G$ on $A$, and the subalgebra $A_\alpha$ of $B$,
defined above,
consisting of the functions $x \mapsto \alpha_x(a)$. Here $X = G$,
and we set ${\Cal M}_\a^+ = {\Cal M}_G^+\cap A_\a$, and similarly
for ${\Cal N}_\a$, ${\Cal M}_\a$, and $\varphi_\a$. But we often
tacitly identify $A$ with $A_\a$.

\definition{1.10 Definition} The action $\alpha$ of $G$ on $A$ is said to be
{\it
integrable} if $A_{\alpha}$ is an integrable subalgebra of $C_b(G, A)$,
that is, if ${\Cal M}_\a$ is dense in $A$.
\enddefinition

We remark that this definition is very close to that given in the sentence
before
$18.20$ of \cite{S} for the setting of a Kac algebra acting on a von~Neumann
algebra.  The case of a group acting on a von~Neumann algebra is then
discussed in
$18.20$ of \cite{S}.

A question which I have not been able to resolve is whether, given an
integrable action $\alpha$ of $G$ on $A$, and given a closed subgroup
$H$ of $G$, the restriction of $\alpha$ to $H$ must always be integrable.
This question is closely related to the notion of ``strongly subgroup
integrable'' introduced in definition 2.17 of \cite{Rf6} and
discussed there.

The following observation about integrable actions is motivated by well-known
considerations in topological dynamics concerning wandering sets.  (See Theorem
$6.15$ of \cite{W}.)

\proclaim{1.11 Proposition} Let $\alpha$ be an action of $G$ on $A$.
Suppose that
$G$ is not compact.  Then every $\alpha$-invariant state on $A$ has value
$0$ on
all of ${\Cal M}_{\alpha}$.  In particular, if $\alpha$ is integrable then
there
are no $\alpha$-invariant states on $A$.
\endproclaim

\demo{Proof} Let $\omega$ be an $\alpha$-invariant state on $A$, and let $a \in
{\Cal M}_{\alpha}$.  By Proposition $1.8$ the function $x \mapsto
\omega(\alpha_x(a)) = \omega(a)$ must be integrable.  Since $G$ is not compact,
$\omega(a) = 0$. \qed
\enddemo

We now return to the case of a general action $\alpha$.  By transport of
structure,
each $\alpha_x$ extends to an automorphism of $A''$, still denoted by
$\alpha_x$,
though the corresponding action of $G$ will usually not be continuous for the
norm.  We will let $(A'')^{\alpha}$ denote the fixed point subalgebra of
$A''$ for
this action.

Now for any $a \in {\Cal M}_{\alpha}^+$ and $x \in G$ the net
$\{\alpha_x(\alpha_f(a))\}_{f \in {\Cal B}}$ must have
$\alpha_x(\varphi(a))$ as
l\.u\.b.  Let $L_x$ denote the usual left translation operator on functions
on $G$
defined by $(L_xf)(y) = f(x^{-1}y)$.  Then for $f \in L^1(G)$ we have
$\alpha_x\alpha_f = \alpha_{L_xf}$.  Furthermore, $L_x$ on ${\Cal B}$ is
clearly an
order automorphism of ${\Cal B}$.  Thus the l\.u\.b\. of
$\{\alpha_x(\alpha_f(a))\}$ must still be $\varphi(a)$.  Consequently
$\alpha_x(\varphi(a)) = \varphi(a)$.  We have thus obtained:

\proclaim{1.12 Proposition} The operator-valued weight $\varphi_{\alpha}$ has
values in $(A'')^{\alpha}$.
\endproclaim

We now examine the action of $\alpha$ on elements of ${\Cal M}_{\alpha}$.  Let
$\Delta$ denote the modular function of $G$, with the convention that
$$
\int f(xy)dx = \Delta(y^{-1}) \int f(x)dx, \hskip 0.7 in
\int f(x^{-1})dx = \int f(x)\Delta(x^{-1})dx.
$$
For $x \in G$ let $R_x$ be the operator of right translation on functions
on $G$
defined by $(R_xf)(y) = f(yx^{-1})$.  We choose this convention both
because then
$R_x$ is an order automorphism of ${\Cal B}$ (not necessarily preserving the
$L^1$-norm), and because a simple calculation shows that for $f \in L^1(G)$
and $a
\in A$ we have
$$
\alpha_f(\alpha_x(a)) = \Delta(x)^{-1}\alpha_{R_xf}(a).
$$
Arguing as we did for Proposition $1.12$ we obtain:

\proclaim{1.13 Proposition} Let $a \in {\Cal M}_{\alpha}$.  Then
$\alpha_x(a) \in
{\Cal M}_{\alpha}$ for any $x \in G$, and
$$
\varphi_{\alpha}(\alpha_x(a)) = \Delta(x)^{-1}\varphi_{\alpha}(a).
$$
\endproclaim

\proclaim{1.14 Corollary} The left ideal ${\Cal N}_{\alpha}$ is carried
into itself
by $\alpha$.
\endproclaim

\demo{Proof} If $a \in {\Cal N}_{\alpha}$ then $a^*a \in {\Cal M}_{\alpha}$, so
that $\alpha_x(a)^*\alpha_x(a) = \alpha_x(a^*a) \in {\Cal M}_{\alpha}$ for
any $x
\in G$. \qed
\enddemo

It is natural to define an $(A'')^{\alpha}$-valued inner-product on ${\Cal
N}_{\alpha}$ by
$$
\<a,b\>_{\alpha} = \varphi_{\alpha}(a^*b).
$$
Note that ${\Cal N}_{\alpha}$ will not in general be a right module over
$(A'')^{\alpha}$, so that ${\Cal N}_{\alpha}$ need not be a Hilbert
$C^*$-module. One can extend ${\Cal N}_{\alpha}$ to get 
a Hilbert $C^*$-module by passing to a
suitable von~Neumann subalgebra of $A''$, but we will not pursue this matter
here.  In any case, we do have a corresponding norm on ${\Cal N}_{\alpha}$
defined
by $\|\<a,a\>_{\alpha}\|^{1/2}$, where the norm in this expression is that
of $A''$.

\demo{\bf 1.15 Notation} For $x \in G$ we define an operator, $U_x$, on ${\Cal
N}_{\alpha}$ by
$$
U_xa = \Delta(x)^{1/2}\alpha_x(a).
$$
\enddemo

A simple calculation shows that $U_x$ is ``unitary'' in the sense that
$$
\<U_xa,U_xb\>_{\alpha} = \<a,b\>_{\alpha}
$$
for $a,b \in {\Cal N}_{\alpha}$.  We obtain in this way a group
homomorphism {}from
$G$ into the group of ``unitary'' operators on ${\Cal N}_{\alpha}$.  It is not
clear to me how often this homomorphism will be strongly continuous for the
norm
defined above.  This seems to be quite a delicate matter to ascertain in
various
examples.  This is closely related to:

\demo{\bf 1.16 Question} Under what circumstances will it be true that for
every
finite measure $\mu$ of compact support on $G$ we have
$
\alpha_{\mu}(a) \in {\Cal N}_{\alpha}$ if  $a 
\in {\Cal N}_{\alpha}
$
(where $\alpha_{\mu}$ is the integrated form of $\alpha$)?
\enddemo

We can show that this is true for $a \in {\Cal M}_{\alpha}$, but we will not
pursue this matter here.

For use in the next section we now examine to some extent what
integrability means
in the commutative case.

\proclaim{1.17 Proposition} Let $\alpha$ be an action of $G$ on the locally
compact space $M$, and so on the $C^*$-algebra $C_{\infty}(M)$ of functions
vanishing at infinity.  If $\alpha$ on
$C_{\infty}(M)$ is
integrable, then every $\alpha$-orbit in $G$ is closed, and the stability
subgroup
of each point of $M$ is compact.
\endproclaim

\demo{Proof} Suppose the $\alpha$-orbit of $m_0 \in M$ is not closed, so
that it
has a limit point $n$ which is not in the orbit.  Choose $f \in C_c(M)^+$ 
(functions of compact support) such
that $f(n) > 1$.  Let $U = \{m \in M: f(m) > 1\}$, so that $U$ is a
neighborhood
of $n$.  By the joint continuity of the action, we can find a symmetric open
precompact neighborhood ${\Cal O}$ of $e_G$ (the identity element of $G$)
and a neighborhood $V$ of $n$ such
that $\alpha_{\Cal O}(V) \subset U$.  Choose a sequence $\{x_j\}$ in $G$ by
induction as follows.  Set $x_1 = e_G$ (the identity element).  
If $x_1,\dots,x_k$ have been
chosen, let
$C_k$ be the closure of the union of the $x_j^{-1}{\Cal O}^2$ for $j \le
k$.  Then
$C_k$ is a compact set.  Thus $\alpha_{C_k}(m_0)$ is a compact subset of
the orbit
of $m_0$, and so can not have $n$ in its closure.  Thus we can find
$x_{k+1} \in
G$ such that $x_{k+1}^{-1} \notin C_k$ and $\alpha_{x_{k+1}}^{-1}(m_0) \in V$.

Notice that since ${\Cal O}$ is symmetric, all the sets $x_j^{-1}{\Cal O}$ are
disjoint.  But if $y \in {\Cal O}$, then
$$
\alpha_{(x_jy)}^{-1}(m_0) = \alpha_y^{-1}\alpha_{x_j}^{-1}(m_0) \in U,
$$
so that
$$
(\alpha_{x_jy}f)(m_0) \ge 1.
$$
That is, the function $x \mapsto (\alpha_x(f))(m_0)$ is non-negative, and has
value $\ge 1$ on each of the disjoint sets $x_n{\Cal O}$, which all have
the same
non-zero Haar measure.  Thus this function can not be integrable.  If we view
evaluation at $m_0$ as a continuous linear functional on $A$, then {}from
Proposition $1.8$ it follows that $f \notin {\Cal M}_{\alpha}^+$.  But ${\Cal
M}_{\alpha}$ is a hereditary $*$-subalgebra, which in the commutative case
means
an ideal.  Thus if ${\Cal M}_{\alpha}$ were dense it would have to contain
$C_c(M)$.  Thus $\alpha$ is not integrable.

We now show that stability groups are compact.  Let $m \in M$, and pick $f \in
C_c(M)^+$ such that $f(m) > 1$.  Let $g(x) = f(\alpha_x^{-1}(m))$.  Then
there is
a compact symmetric neighborhood ${\Cal O}$ of $e_G$ on which $g \ge 1$.  Let
$G_m$ denote the stability subgroup of $m$.  Then $g(xs) = g(x)$ for $x \in
G$, $s
\in G_m$.  If $G_m$ is not compact, it is easily seen {}from this that $g$
can not
be integrable. \qed
\enddemo

I have not noticed simple conditions which are simultaneously necessary and
sufficient for
$\alpha$ to be integrable.  It is not sufficient just to have the orbits be
closed
and the stability subgroups be compact.  This is seen by the following example,
which is a slightly more complicated version of the example at the very end of
Philip Green's article \cite{G}.  We make our example yet slightly more
complicated than needed here so that we can also use it in the next section to
illustrate a point there.

\example{1.18 Example} The space $M$ is a closed subset of ${\Bbb R}^3$,
and the
group $G$ is ${\Bbb R}$.  The action is free, with all orbits closed.  The
orbit
space $M/\alpha$ is a compact Hausdorff space consisting of a countable
number of
points, which is discrete except for one limit point.  This limit orbit is the
``$y$-axis'' $\{(0,s,0): s \in {\Bbb R}\}$, with the action $\alpha$ of
${\Bbb R}$
on it being by translation.  We denote this orbit by ${\Cal O}_*$, and we
let $p_*
= (0,0,0)$, which is one of its points.  We label the other orbits by strictly
positive integers, $n$, and we denote the $n$-th orbit by ${\Cal O}_n$.
Part of
the data specifying these orbits consists of a strictly decreasing sequence
$\{b_n\}$ of real numbers which converges to $0$.  Let $p_n = (b_n,0,0)$.  Then
$p_n$ will be in ${\Cal O}_n$.  Up to equivariant homeomorphism the example
will be
independent of the choice of $\{b_n\}$.  However, it does depend on the
next piece
of data, which is an assignment to each $n$ of a strictly positive integer,
$L_n$,
which should be thought of as a ``repetition number''.  However, the
example will be independent of the choise of the
next piece of data, which is an assignment to each $n$ of a
strictly decreasing finite sequence $\{c_j^n\}$ of length $L_n$, with
$b_{n+1} <
c_j^n < b_n$.  Let $q_j^n = (c_j^n,0,0)$.  Each of the points $q_j^n$, $j =
1,\dots,L_n$, will be in the orbit ${\Cal O}_n$.  At one place below it is
convenient to set $c_0^n = b_n$.  We specify ${\Cal O}_n$ and the action
$\alpha$
by saying where $\alpha$ takes $p_n$.  For $t \in {\Bbb R}$ we set
$$
\alpha_t(p_n) = \cases
(b_n,t,0) &t \in (-\infty,n]. \\
(c_{L_n}^n,s,0) &s \in (-n,+\infty),\ t = s + L_n(2n+1). \\
(c_j^n,s,0) &s \in (-n,n],\ t = s + j(2n+1), \\
&1 \le j < L_n, \text{ void if } L_n = 1. \\
((1-s)c_j^n + sc_{j+1}^n,n \cos(\pi s),n \sin(\pi s)) &s \in (0,1],\ t = s
+ n +
j(2n+1), \\
&0 \le j \le L_n - 1,\ c_0^n = b_n.
\endcases
$$

If one contemplates the facts that whenever points of $m$ are in the
$x$-$y$-plane
and not about to leave it they move parallel to the $y$-axis with unit
speed, and
that, as $n$ increases, the $y$-coordinates where points enter and leave the
$x$-$y$-plane go to $-\infty$ and $+\infty$, 
one sees that this action is indeed
jointly continuous, and that the properties stated at the beginning are
satisfied.
In particular, for any $f \in C_c(M)$, the support of $f$ will meet any
given orbit
in a compact set, and thus
$$
\int_G f(\alpha_x(m))dx
$$
will be finite for each $m \in M$.

Choose now an $f \in C_c({\Cal O}_*)^+$, supported strictly inside $\{0\}
\times
[-1/2,1/2] \times \{0\}$, and such that its integral over ${\Cal O}_*$, i\.e\.
$\int f(\alpha_{-t}(p_*))dt$, $= 1$.  
Extend $f$ to a function in $C_c({\Bbb R}^2 \times
\{0\})^+$, still denoted by $f$, in such a way that the support of $f$ is
contained
in the disk of radius $1/2$ about the origin, and that for some
$\varepsilon > 0$
this extended $f$ is independent of the $x$-coordinate.  We restrict $f$ to
$M$ and
still denote it by $f$.  Then as soon as $n$ is large enough that $b_n <
\varepsilon$, the restriction of $f$ to ${\Cal O}_n$ will, on $G$, look
like $L_n +
1$ copies of $f$ with disjoint supports.  Consequently, for each large $n$
we have
$$
\int f(\alpha_{-t}(p_n))dt = L_n + 1.
$$
In other words, for large enough $n$,
$$
\int \alpha_t(f)dt = \cases
1 &\text{on ${\Cal O}_*$} \\
L_n + 1 &\text{on ${\Cal O}_n$.}
\endcases
$$
In particular, if the sequence $\{L_n\}$ is not bounded, then $f$ is not
integrable, so that ${\Cal M}_{\alpha}$ is not dense, and so the action
$\alpha$ is
not integrable.  On the other hand, if the sequence $\{L_n\}$ is bounded,
then one
can check that $\alpha$ is integrable.

\endexample

We remark that by examining the foliations of the plane (which come {}from 
actions of ${\Bbb R}$), as studied in \cite{Wn}, we obtain an abundance of 
examples of integrable actions $\alpha$ on $C_{\infty}(M)$ such that
$M/\alpha$ is not Hausdorff (but the actions are free, with closed orbits).
Thus, integrability does not imply that $M/\alpha$ is Hausdorff.

Anyway, we are left with:

\demo{\bf 1.19 Question} What are conditions just in terms of an action
$\alpha$ of a group $G$
on a space $M$ which are simultaneously necessary and sufficient 
for $\alpha$ on $C_{\infty}(M)$ to be integrable?
\enddemo

We conclude this section with an important property of integrable actions with
respect to tensor products.  Many special cases of this property are
employed in
the literature.  (See e\.g\. $18.21$ of \cite{S}.)

\proclaim{1.20 Proposition} Let $\alpha$ and $\beta$ be actions of $G$ on the
$C^*$-algebras $A$ and $B$, and let $\alpha \otimes \beta$ denote the
corresponding
action on their maximal, or minimal, tensor product, $A \otimes B$.  If
$\alpha$ is
integrable, then so is $\alpha \otimes \beta$.
\endproclaim

\demo{Proof} It does not hurt to adjoin a unit to $B$ if it does not have
one.  So
we assume that $B$ is unital.  Let $a \in {\Cal M}_{\alpha}^+$.  Then
$$
(\alpha \otimes \beta)_\lambda(a \otimes 1_B) = \alpha_\lambda(a) \otimes 1 \le
k_a\|\lambda\|_{\infty}
$$
with our earlier notation.  Thus $a \otimes 1_B \in {\Cal
M}_{\alpha\otimes\beta}^+$.  But for any $b \in B^+$ we have $a \otimes b \le
\|b\|(a \otimes 1_B)$.  Since ${\Cal M}_{\alpha\otimes\beta}^+$ is
hereditary, it
follows that $a \otimes b \in {\Cal M}_{\alpha\otimes\beta}^+$.  Since ${\Cal
M}_{\alpha}$ is dense in $A$ by assumption, it follows that ${\Cal
M}_{\alpha\otimes\beta}$ is dense in $A \otimes B$. \qed
\enddemo

\section{2.  Strict C*-weights}

We recall that if $\alpha$ is an action of $G$ on a locally compact space
$M$, then
$\alpha$ is said to be {\it proper} if the map $(m,x) \to (m,\alpha_x(m))$
{}from $M
\times G$ to $M \times M$ is proper, in the sense that preimages of compact
sets
are compact.  It is well known \cite{Bo} that in this case the orbit space,
$M/\alpha$, with the quotient topology {}from $M$, is locally compact
Hausdorff.
The functions in $C_{\infty}(M/\alpha)$ can be viewed as 
functions in $C_b(M)$
which are $\alpha$-invariant.  Here $C_b(M)$ is the algebra of bounded
continuous
functions on $M$, and it is the multiplier algebra of $C_{\infty}(M)$.  It is
well-known ($2.4$ of \cite{Pn}) that if $h \in C_c(M)$, and if we set
$$
\psi(h)(m) = \int h(\alpha_x^{-1}(m))dx
$$
for every $m \in M$, then $\psi(h)$ is a function in $C_{\infty}(M/\alpha)
\subseteq
C_b(M)$.  It is natural to write
$$
\psi(h) = \int \alpha_x(h)dx,
$$
but as before, the integrand is not integrable in the usual sense if $G$ is not
compact.  But if we consider $\alpha_\lambda(h)$ 
for $\lambda \in {\Cal B}$ as in the
previous
section, it is easily seen that $\alpha_\lambda(h)$ 
converges to $\psi(h)$ in the
strict
topology, that is, $k\alpha_\lambda(h)$ 
converges to $k\psi(h)$ in (uniform) norm for
every $k \in C_{\infty}(M)$.  See the discussion early in section~1 of
\cite{Rf7}.
For the definition and basic properties of the strict topology 
see \cite{La,Pe2}.

There are a number of situations in which an action of a group on a
non-commutative
$C^*$-algebra seems to be proper in some sense.  I tried to give an appropriate
definition in \cite{Rf7}.  The definition given there was adequate to treat
some
interesting examples, but it assumed the existence of a dense subalgebra with
certain properties, and so was not intrinsic.  I propose here a tentative
intrinsic
definition, which is essentially one almost explicitly suggested by Exel in
section~6 of \cite{E1}.  The difference is that here we emphasize the order
properties, paralleling the development in our first section, while Exel
emphasizes
unconditional integrability.

As suggested by the above discussion, this matter leads to weights with
values in
$C^*$-algebras.  We mentioned in the previous section that such weights have
recently
been introduced by Kustermans \cite{Ku} for use in connection with quantum
groups
(though so far I have not seen how to use his ``regularity'' condition in the
present context).   Much as we need here, he treats weights on a
$C^*$-algebra $B$
with values in $M(A)$ for another $C^*$-algebra $A$.  (Here $M(A)$ denotes the
multiplier algebra \cite{Pe2} of $A$.) Our basic context is as follows:

\definition{2.1 Definition} Let $B$ and $A$ be $C^*$-algebras, and 
let $P = \{p_{\lambda}\}$ be an
increasing net of positive operators {}from $B$ into $M(A)$.
We say that $b \in B^+$ is $P$-{\it proper} if
the net
$\{p_{\lambda}(b)\}$ converges in the strict topology to an element,
$\psi_P(b)$,
of $M(A)$.  We denote the set of $P$-proper elements of $B^+$ by ${\Cal
P}_P^+$.
\enddefinition

It is clear that ${\Cal P}_P^+$ is a cone, and that $\psi_P$ is ``linear'' on
${\Cal P}_P^+$.

For use in dealing with this definition we now recall several 
of the basic facts about the strict topology which Kustermans
obtains, in a form suitable for our needs here. The considerations here
parallel somewhat the strong and weak operator topologies. A small novelty 
is our explicit definition of the ``weak-strict'' topology. (It has
been used implicitly in earlier work.)

\definition{2.2 Definition} We say that a net $\{m_{\lambda}\}$ in
$M(A)$
converges in
the {\it weak-strict topology} to $n \in M(A)$ if the net $\{am_{\lambda}c\}$
converges in norm to $anc$ for every $a,c \in A$.  By polarization it
suffices to
consider
$\{a^*m_{\lambda}a\}$ and $a^*na$. We will say that a net $\{m_{\lambda}\}$ is
{\it weak-strict Cauchy} if for every $a ,c \in A$ the net
$\{am_{\lambda}c\}$ in
$A$ is norm-Cauchy.  Again, it suffices to examine the nets
$\{a^*m_{\lambda}a\}$ for $a \in A$.  
\enddefinition

It is clear that if a net is strictly convergent, or is strictly Cauchy,
then it is also
so weak-strictly.

\proclaim{2.3 Proposition} Let $\{m_{\lambda}\}$ be an increasing net in
$M(A)^+$
which converges weak-strictly to $n \in M(A)$.  Then $m_{\lambda} \le n$
for all
$\lambda$.  In particular, $\{m_{\lambda}\}$ is bounded in norm.
\endproclaim

\demo{Proof} Fix any $\lambda_0$, and set $k = n - m_{\lambda_0}$ and
$k_{\lambda}
= m_{\lambda} - m_{\lambda_0}$.  Then the net $\{k_{\lambda}\}$ is eventually
positive and converges weak-strictly to $k$.  Thus for any $a \in A$ the net
$\{a^*k_{\lambda}a\}$ is eventually positive and converges in norm to $a^*ka$.
Thus $a^*ka$ is positive for all $a \in A$.  Then it is not hard to see
that $k$ is
positive.  (See Lemma $4.1$ of \cite{La}.) \qed
\enddemo

In lemma $9.3$ of \cite{Ku} Kustermans uses the uniform boundedness principle
several times to show that in the above proposition it suffices to assume that
$\{m_{\lambda}\}$ is weak-strict Cauchy. This observation can be useful in
connection with the following proposition.

\proclaim{2.4 Proposition} (See $9.4$ and $9.5$ of \cite{Ku}.)  Let
$\{m_{\lambda}\}$ be an increasing net in $M(A)^+$ which is weak-strict
Cauchy and
is bounded in norm.  Then $\{m_{\lambda}\}$ is strictly Cauchy, and so
converges
strictly (and so weak-strictly) to an element of $M(A)^+$.
\endproclaim

\demo{Proof} Let $K$ be a bound on $\{m_{\lambda}\}$.  Then for $a \in A$ and
$\lambda > \mu$ we have
$$
\align
\|(m_{\lambda}-m_{\mu})a\|^2 &\le \|(m_{\lambda}-m_{\mu})^{1/2}\|^2
\|(m_{\lambda}
- m_{\mu})^{1/2}a\|^2 \\
&\le K\|a^*(m_{\lambda}-m_{\mu})a\|.
\endalign
$$
Thus $\{m_{\lambda}a\}$ is norm-Cauchy.  By taking adjoints
we see that $\{am_{\lambda}\}$ too is
norm-Cauchy.  Thus $\{m_{\lambda}\}$ is strictly Cauchy, and so converges to a
positive element of $M(A)$, since $M(A)$ is strictly complete \cite{La}. \qed
\enddemo

\proclaim{2.5 Proposition} (Lemma $9.4$ of \cite{Ku}.)  Let
$\{m_{\lambda}\}$ be a
net of elements of $M(A)^+$, and let $m \in M(A)^+$ be such that
$m_{\lambda} \le
m$ for all $\lambda$.  If $\{m_{\lambda}\}$ converges weak-strictly to $m$,
then it
does so strictly.
\endproclaim

\demo{Proof} For any $a \in A$ we have, as above,
$$
\align
\|(m-m_{\lambda})a\|^2 &\le
\|(m-m_{\lambda})^{1/2}\|^2\|(m-m_{\lambda})^{1/2}a\|^2
\\
&\le \|m\| \|a^*(m-m_{\lambda})a\|,
\endalign
$$
For $a$ on the other side, take adjoints. \qed
\enddemo

There is a useful alternate characterization of the $P$-proper elements
in terms of linear functionals on $A$. It is related to the definition
of $\hat \alpha$-integrable elements given on page 269 of \cite{Pe2},
which originated in \cite{OP1, OP2}.
But we call attention to the note at the end of \cite{OP2} which points
out that the definition in \cite{Pe2} is too strong, since it should only
consider the dual, $B'$, of $B$ (in the notation of \cite{Pe2}), not
of $M(B)$. For some later variants see \cite{QR} and its references. We
will use the fact that each element of $A'$ has a canonical extension
to $M(A)$, obtained by viewing $A'$ as the predual of $A''$ and $M(A)$
as canonically embedded in $A''$ (proposition 3.12.3 of \cite{Pe2}).

\proclaim{2.6 Theorem} Let $P = \{p_\lambda\}$ be an increasing net of
positive operators {}from $B$ to $M(A)$, and let $b \in B^+$. Then $b$ is
$P$-proper if and only if there is an $m \in M(A)^+$ such that for
every positive linear functional, $\omega$, on $A$, viewed also as a
positive linear functional on $M(A)$, the net $\{\omega(p_\lambda(b))\}$
converges to $\omega(m)$.
\endproclaim

\demo{Proof} Suppose that $b$ is $P$-proper. For any bounded linear
functional $\omega$ on $A$ and any bounded net $\{n_\lambda\}$ in $M(A)$
which converges strictly to $m \in M(A)$ the net $\{\omega(n_\lambda)\}$
converges to $\omega(m)$. This follow {}from the fact that $\omega$
can be approximated in norm by linear functionals of the form
$\omega_u$ defined by $\omega_u(a) = \omega(u^*au)$ for $u \in A$.
(See e.g. the proof of theorem 3.12.9 of \cite{Pe2}.) {}From this it follows
that if $b$ is $P$-proper, then $\{\omega(p_\lambda(b)\}$ converges
to $\omega(\psi_P(b))$.

Suppose conversely that there is an $m \in M(A)^+$ as in the statement
of the theorem. For any positive $\omega$ and any $c \in A$ let $\omega_c$
be defined by $\omega_c(a) = \omega(c^*ac)$. Since $\omega_c$ is positive,
$\omega_c(p_\lambda(b))$ converges by hypothesis to $\omega_c(m)$, that
is, $\omega(c^*p_\lambda(b)c)$ converges to $\omega(c^*mc)$. But now
$c^*mc \in A$. Let $Q(A)$ denote the quasi-state space \cite{Pe2} 
of $A$, consisting of those positive $\omega \in A'$ such that
$\|\omega\| \leq 1$. Note that $Q(A)$ is weak-* compact. For any
$d \in A$ let $\hat d$ denote $d$ viewed as a function on $Q(A)$, so
that $\hat d$ is continuous (and affine). With this notation, 
$(c^*p_\lambda(b)c)\ \hat{}\ \ $ is an increasing net of continuous positive
functions on $Q(A)$, which as we saw above converges pointwise
to the continuous function $(c^*mc)\ \hat{}\ \ $. By Dini's theorem it 
follows that the convergence is uniform. But $Q(A)$ determines the
norm of elements of $A^+$. It follows that the net $\{c^*p_\lambda(b)c\}$
converges in norm to $c^*mc$, that is, $\{p_\lambda(b)\}$ converges
weak-strictly to $m$. {}From Proposition 2.5 it follows 
that $\{p_\lambda(b)\}$ converges to $m$ strictly. Thus $b$ is $P$-proper
as desired.     \qed
\enddemo

The following lemma is motivated by, and very closely related to,
proposition $6.6$
of \cite{E1}and to the comments in $7.8.4$ of \cite{Pe2} and $2.4$ of
\cite{OP1}.
See also lemma 3.5 of \cite{Qg} and lemma $4.1$ of \cite{Ku}.

\proclaim{2.7 Key Lemma} Let $P = \{p_{\lambda}\}$ be an increasing net of
positive operators {}from $B$ to $M(A)$.  Then the cone ${\Cal P}_P^+$ of
$P$-proper elements of $B^+$ is hereditary.
\endproclaim

\demo{Proof} Let $b \in {\Cal P}_P^+$, and let $b_0 \in B$ with $0 \le b_0
\le b$.
The net $\{p_{\lambda}(b)\}$ is bounded above by Proposition $2.3$, and so
the net
$\{p_{\lambda}(b_0)\}$ is bounded above.  Since the latter net is
increasing, to
show that $b_0 \in {\Cal P}_P^+$ it suffices by Proposition $2.4$ to show that
$\{p_{\lambda}(b_0)\}$ is weak-strict Cauchy.  Now for any $a \in A$ and any
$\lambda > \mu$
$$
a^*(p_{\lambda}(b_0) - p_{\mu}(b_0))a = a^*((p_{\lambda}-p_{\mu})(b_0))a
\le a^*((p_{\lambda}-p_{\mu})(b))a.
$$
But $\{p_{\lambda}(b)\}$ converges strictly, and so is weak-strict Cauchy. \qed
\enddemo

According to the properties of hereditary cones given in section~1, if we set
$$
{\Cal Q}_P = \{b \in B: b^*b \in {\Cal P}_P\},
$$
then ${\Cal Q}_P$ is a left ideal in $B$.  Let ${\Cal P}_P$ denote the
linear span
of ${\Cal P}_P^+$.  Then ${\Cal P}_P = {\Cal Q}_P^*{\Cal Q}_P$, ${\Cal P}_P
\cap
B^+ = {\Cal P}_P^+$, and $\psi_P$ extends uniquely to a positive 
linear map {}from ${\Cal P}_P$ into $M(B)$.

It is clear that every $P$-proper positive element is $P$-bounded.  Thus ${\Cal
P}_P^+ \subseteq {\Cal M}_P$ in the notation of the previous section.
Furthermore, since
any non-degenerate representation of $B$ extends uniquely to one of $M(B)$,
under
which a strictly convergent net is strong operator convergent, $\psi_P$
will be the
restriction of the weight $\varphi_P$ to the $P$-proper elements.  The above
considerations suggest:

\definition{2.8 Definition} Let $A$ and $B$ be $C^*$-algebras.  By a
$C^*$-valued
weight on $B$ {\it towards} $A$ we mean a $C^*$-weight $\psi$ on $B$
(Definition
$1.2$) with values in $M(A)$.  Let ${\Cal P}^+$ denote the domain of
$\psi$.  We say that $\psi$ is {\it lower semi-continuous} if
whenever $b \in {\Cal P}^+$ and 
$\{b_{\mu}\}$ is an increasing net in $B^+$ which converges in norm to $b$,
then
the net $\{\psi(b_{\mu})\}$ converges strictly to $\psi(b)$.
We
say that $\psi$ is {\it strict} if there is an increasing net
$\{p_{\lambda}\}$ of
bounded positive maps {}from $B$ to $M(A)$ for which

\roster
\item"1)" ${\Cal P}^+ = \{b \in B^+: \text{ the net $\{p_{\lambda}(b)\}$ is
strictly Cauchy}\}$.
\item"2)" $\psi(b) = \text{ strict-lim} \{p_{\lambda}(b)\} \text{ for }
b \in
{\Cal P}^+$.
\endroster

If $B = A$ we will just say that $\psi$ is a {\it lower semi-continuous},
or {\it strict, $C^*$-valued
weight} on
$A$.  If $\psi(b) = 0$ only when $b = 0$, we say that $\psi$ is {\it faithful}.
\enddefinition

We remark that we do not require that ${\Cal P}$ be dense, unlike
definition $3.2$
of \cite{Ku}. Nor do we require {\it complete} positivity. (We will
assume it explicitly when we need it.).  
The exact relationship between our definition of ``strict''
weights and Kustermans' definition of lower semi-continuous weights in
definition 3.2 of \cite{Ku} remains to be worked out. We note that
in \cite{Ku} each $p_{\lambda}$ is required to be
``strict''
as defined in \cite{La}.  This has some technical advantages, 
but in Proposition
$4.11$ we will see that we have an even stronger property in our
group-action case.

\proclaim{2.9 Proposition} (Basically $3.5$ of \cite{Ku}.)  Any
strict
$C^*$-valued weight {}from $B$ toward $A$ is lower semi-continuous.
\endproclaim

\demo{Proof} Let $b$ and $\{b_\mu\}$ be as in the definition above of
lower semi-continuity. By Proposition $2.5$ it suffices to show weak-strict
convergence.  Let
$a \in A$, and let $\varepsilon > 0$ be given.  Choose $\lambda$ so that
$$
\|a\psi(b)a^* - ap_{\lambda}(b)a^*\| < \varepsilon/2.
$$
Choose $\mu_0$ such that if $\mu > \mu_0$ then
$$
\|ap_{\lambda}(b)a^* - ap_{\lambda}(b_{\mu})a^*\| < \varepsilon/2.
$$
Since we have
$$
ap_{\lambda}(b_{\mu})a^* \le a\psi(b_{\mu})a^* \le a\psi(b)a^*,
$$
it follows that for $\mu > \mu_0$ we have
$$
\|a\psi(b_{\mu})a^* - a\psi(b)a^*\| < \varepsilon.
$$
\qed
\enddemo

We remark that in the situation described above we can not expect that
$\psi(b_\lambda)$ will converge to $\psi(b)$ in norm. A simple example,
which we will use again later, goes as follows:

\example{2.10 Example} Let $G = {\Bbb Z}$ act by translation, $\tau$, on 
itself, and so on $C_\infty({\Bbb Z})$. For each $n \geq 1$ choose
$f_n \in C_c^+({\Bbb Z})$ such that $\|f_n\|_\infty \leq 1/n$ but
$\sum_k \tau_k(f_n) = 1$ strictly in $M(C_\infty({\Bbb Z})$. Let $\alpha$
be the (proper) action of ${\Bbb Z}$ by translation in the first
variable on ${\Bbb Z} \times {\Bbb N}$, and 
so on $B = C_\infty({\Bbb Z} \times {\Bbb N})$. Let $P = \{p_\lambda\}$
come {}from $\alpha$ as at the beginning of Section 1, with
corresponding $\psi$. Let $g \in B$ be defined by $g(m,n) = f_n(m)$.
It is easily seen that $g$ is $P$-proper, and that $\psi(g) = 1$. Define
$g_j$ to agree with $g$ for $n \leq j$ and have value $0$ otherwise.
Then the increasing sequence 
$\{g_j\}$ converges to $g$ in norm, while $\{\psi(g_j)\}$
converges to $\psi(g)$ strictly, but not in norm.
\endexample

\section{3. The case of $C_b(X,A)$}

Throughout this section we let $B=C_b(X,A)$, and we assume
that we have a positive Radon measure specified on $X$, in terms of
which the $p_{\l}$'s are defined as near the beginning of Section 1.
We now denote the hereditary cone of $P$-proper elements by ${\Cal P}_X^+$,
with corresponding $\psi_X$, ${\Cal P}_X$  and ${\Cal Q}_X$.  
We consider here some aspects which are
special to this situation.

As discussed in the previous section, ${\Cal Q}_X$ will always be a
left ideal in $B$. We now consider action on the right.
Let $f\in\Cal Q_X$ and $m\in M(A)$, so
that $fm\in B$. For $\l\in\Cal B$ we have
$$
p_{\l}((fm)^*fm)=\int m^*f(x)^*f(x)m\lambda (x)dx =
m^*p_{\l}(f^*f)m \  ,
$$
and we know that $a^*m^*p_{\l}(f^*f)ma$ converges up in norm to
$a^*m^*\psi(f^*f)ma$. Thus $p_{\l}((fm)^*(fm))$ converges weak-strictly
to $m^*\psi(f^*f)m$, and so converges strictly by Proposition 2.5.
By definition it follows that $fm\in \Cal Q_X$, so that $\Cal Q_X$ 
is a right $M(A)$-module, and
$$
\psi_X(m^*f^*fm)=m^*\psi_X(f^*f)m   .
$$
Let also $g\in\Cal Q_X$.
Then $g^*fm\in\Cal P_X$, and so $\psi(g^*fm)$ is the strict limit of
$\{p_{\l}(g^*fm)\}$. Consequently for $a\in A$ we have the norm limits
$$
a\psi_X(g^*fm) =\lim ap_{\l}(g^*fm) =\lim (ap_{\l}(g^*f))m
=a\psi_X(g^*f)m \ .
$$
Thus $\psi_X(g^*fm)=\psi_X(g^*f)m$. Finally, since every element of
$\Cal P_X$ is a finite linear combination of elements of form
$g^*f$ for $f,g\in \Cal Q_X$, we see that we have obtained:

\proclaim{3.1 \ Proposition} Both $\Cal Q_X$ and $\Cal P_X$ are right
$M(A)$-modules, and
$$
\psi_X(fm)=\psi_X(f)m
$$
for $f\in\Cal P_X$ and  $m\in M(A)$.
\endproclaim

We can now define an $M(A)$-valued inner-product on $\Cal Q_X$
by
$$
\< f, \ g\>_X =\psi(f^*g) \ .
$$
This evidently makes $\Cal Q_X$ into a right $C^*$-module over $M(A)$.
(See, e.g.~\cite{La} for the definition.) Consequently we have the
following version of the Cauchy-Schwarz inequality (proposition 2.9 of
\cite{Rf2}, or proposition 1.1 of \cite{La}):

\proclaim{3.2 \ Proposition} For $f,g\in\Cal Q_X$ we have
$$
\psi(f^*g)^*\psi(f^*g) \leq \|\psi(f^*f)\|\psi(g^*g)
$$
in $M(A)$. Consequently the expression $\|f\|_X=\|\<f,f\>_X\|^{1/2}$
defines a norm on $\Cal Q_X$.
\endproclaim

We will later have use for the following technical consequence.

\proclaim{3.3 \ Proposition} With notation as above, for any
$\l\in\Cal B$ and any $f, \ g\in\Cal Q_X$ we have
$$
\|p_{\l}(f^*g)\| \leq 4\|\<f, \ f\>_X\| \ \|\<g, \ g\>_X\| \ .
$$
\endproclaim

{\smc Proof.} By polarization
$$
\align
4\|p_{\l}(f^*g)\| & \leq \sum^3_{k=0}
\|p_{\l}((f + i^kg)^*(f+i^kg)\| \\
& \leq \sum^3_{k=0} \|\psi((f+i^kg)^*(f+i^kg))\|
\leq 4(\|f\|_X+ \|g\|_X)^2 \ ,
\endalign
$$
where for the last inequality we have used the last part of Proposition 3.2.
Now replace $f$ by $f/\|f\|_X$ and $g$ by $g/\|g\|_X$ to obtain
the desired inequality. \qed

\bigskip
We deduce next the application of Theorem 2.6 to the present case.

\proclaim{3.4 \ Theorem} Let $f\in B^+$. Then $f\in\Cal P_X$ if and only
if there is an $m\in M(A)^+$ such that for every positive
$\ \omega\in A'$ the function $x\mapsto\omega(f(x))$ is integrable in the
ordinary sense and
$$
\int\omega(f(x))dx =\omega(m) \ .
$$
\endproclaim

{\smc Proof.} If $f\in\Cal P_X$
then by Proposition 1.8 the function $x\mapsto\omega(f(x))$ is
integrable and the above equation holds by the comment just before
Definition 2.8.

Conversely, if $m$ exists as in the statement of the proposition, then
for any $\lambda\in\Cal B$ we have $\lambda(x)\omega(f(x))\leq \omega(f(x))$.
Since $x\mapsto\omega(f(x))$ is integrable and the $\l$'s converge up
to 1 pointwise, we can pass to a suitable sequence of $\l$'s to which
we can apply the monotone convergence theorem to conclude that
the net $\{\omega(p_{\l}(f))\}$ converges up to $\omega(m)$.
We are now exactly in position to apply Theorem 2.6. \qed

\section{4. Proper Actions}

We now return to the case of an action $\alpha$ of a group $G$ on
a C*-algebra $A$. As suggested earlier, we can view $A$ as embedded
in $C_b(G,A)$ by sending $a \in A$ to the function $x \mapsto \alpha_x(a)$.
We can then apply our earlier discussion to this subalgebra. However,
for our later discussion of morphisms we need a slight generalization
of this situation. The action $\a$ extends to an action, still denoted
by $\a$, on $M(A)$, which need not be strongly continuous.
If we view $M(A)$ as included in $A''$ \cite{Pe2}, this action
is just
the restriction of that on $A''$ used in section 1.  
We let $M(A)_e$ denote the ``$\a$-essential'' part of $M(A)$ for this
action, that is, the C*-subalgebra of elements on which $\a$
is strongly continuous. Then we can extend $p_\lambda$ to $M(A)_e$
by the same formula as before.

\definition{4.1 Definition} We say that $n \in (M(A)_e)^+$ 
is $\alpha$-{\it proper}
if there exists an $m \in M(A)$ such that the 
net $\{p_\lambda(n)\}_{\lambda \in {\Cal B}}$ converges strictly to $m$,
where now 
  $$p_\lambda(n) = \int_G \lambda(x)\alpha_x(n) dx.$$
We denote the hereditary cone of $\alpha$-proper positive 
elements by ${\tilde {\Cal P}}_\alpha^+$ (notice that Lemma 2.7
applies here), and the corresponding strict
C*-weight {}from $M(A)$ towards itself by ${\tilde \psi}_\alpha$. 
We let ${\Cal P}^+ = {\tilde {\Cal P}}^+ \cap A$, a hereditary
cone in $A$.  We have the 
corresponding left ideals ${\tilde {\Cal Q}_\a}$ and ${\Cal Q}_\alpha$, 
and $\ast$-subalgebra ${\tilde {\Cal P}_\a}$ and ${\Cal P}_\alpha$,  
and we let $\psi_\a$ be the restriction 
of ${\tilde \psi}_\a$ to ${\Cal P}_\a$.
\enddefinition

We remark that, in contrast to \cite{Qg, QR}, our restriction to the
$\a$-essential part of $M(A)$ is required in order for $p_\lambda$ to
be defined by an ordinary integral. For example, when $A = C_\infty(\Bbb R)$
and $\a$ is the action of $\Bbb R$ by translation, $M(A)_e$ consists
of the uniformly continuous bounded functions. But there exist functions
in $C_b(\Bbb R) \cap L_1(\Bbb R)$ which are not uniformly continuous.
For such a function $g$ our definition of $p_\lambda(g)$ will not make
sense as it stands. But $x \mapsto \omega(\a_x(g))$ will be
integrable for every finite measure $\omega$. One can develop a more
complicated definition of $p_\lambda$ to handle this kind of situation,
but so far I have not seen a need for this.

We now give an example to show that even for the case of a proper action of
$G$ on
a space $M$, there can be many positive $\alpha$-integrable elements 
in $C_\infty(M)$ which
are not
$\alpha$-proper.

\example{4.2 Example} Let $M = {\Bbb R}$, let $G = 4{\Bbb Z}$ and let
$\alpha$ be
the action of $G$ on $M$ by translation.  This is a proper action.  Let
$g_n$ be
the function on $[0,1]$ defined by
$$
g_n(t) = \cases
nt &\text{for $0 \le t \le 1/n$} \\
1 &\text{otherwise.}
\endcases
$$
Note that the sequence $\{g_n\}$ increases to $\chi_{(0,1]}$, the
characteristic
function of $(0,1]$.  Let $h_1 = g_1$, and for $n \ge 2$ let $h_n = g_n -
g_{n-1}$.  Thus $g_n = \sum_{j=1}^n h_j$.  It is easily seen that
$\|h_n\|_{\infty}
= 1/n$.  Let $k_n$ be the reflection of $h_n$ about $t = 1$, extended to be $0$
outside $[0,2]$.  Then $k_n \in C_c({\Bbb R})$, $\|k_n\|_{\infty} = 1/n$, and
$\sum_{j=1}^n k_j$ converges up to $\chi_{(0,2)}$.  Let $L_n$ be translation by
$4n$.  Set
$$
f = \sum^{\infty} L_nk_n,
$$
where we note that the convergence is uniform.  
Since $\|k_n\|_{\infty} = 1/n$, it follows
that $f
\in C_{\infty}({\Bbb R})$.  Then it is easily seen that $f$ is
$\alpha$-integrable.  But if we identify $M/\alpha$ with the fundamental domain
$[0,4)$, then it is easily seen that
$\sum_{x \in G} \alpha_x(f)$ is $\chi_{(0,2)}$,
which is not in $M(A)$.
\endexample

We remark that a related example appears as example $2.4$ of \cite{FMT}.

We now have the corresponding version of Theorem 3.4. It relates our
situation to definition 3.4 of \cite{Qg} and section 1 of \cite{QR}.

\proclaim{4.3 Theorem} Let $n \in (M(A)_e)^+$.  Then $n \in {\tilde {\Cal
P}}_{\alpha}^+$ if
and only if there is an $m \in M(A)^+$ such that for every positive 
$\omega \in A'$,
viewed
also as a linear functional on $M(A)$, the function $x \mapsto
\omega(\alpha_x(n))$
is integrable on $G$ and
$$
\int \omega(\alpha_x(n))dx = \omega(m).
$$
\endproclaim

\demo{Proof} Suppose that $n \in {\tilde {\Cal P}}_{\alpha}^+$ 
and $\omega \in A'$.  It
follows {}from Proposition $1.8$ and the comments made just before
Definition 2.8 that $x \mapsto \omega(\alpha_x(a))$ is
integrable,
with integral $\omega(\psi(m))$.

Suppose conversely that $x \mapsto \omega(\alpha_x(n))$ is integrable for all
positive $\omega \in A'$, and that there is an $m \in M(A)^+$ such that $\int
\omega(\alpha_x(n))dx = \omega(m)$ for all $\omega$.  Much as in the
proof of
Proposition $1.8$, the net $\{\omega(p_\lambda(n))\}$ converges to $\omega(m)$
as $\lambda$
ranges through ${\Cal B}$.   {}From Proposition $2.6$ it
follows
that $p_\alpha(n)$ converges to $m$ strictly.  Thus $n \in {\tilde {\Cal
P}}_{\alpha}^+$ as
desired. \qed
\enddemo

Exel \cite{E1, E2} defines $a \in A$ to be $\alpha$-integrable if for all $b \in
B$ the
functions $x \mapsto b\alpha_x(a)$ and $x \mapsto \alpha_x(a)b$ are
unconditionally
integrable (meaning that the net $\{\int_E b\alpha_x(a)dx\}$ for $E$
ranging over
precompact subsets of $G$ is norm Cauchy, and similarly for $b$ on the other
side).  The integrals, as $b$ ranges over $A$, then define an element of
$M(A)$.
For general $a \in A$ it is not clear to me whether this implies that $a
\in {\Cal
P}_{\alpha}$.  But for positive elements we have:

\proclaim{4.4 Proposition} Let $a \in A^+$.  Then $a$ is
$\alpha$-integrable in
Exel's sense iff $a \in {\Cal P}_{\alpha}^+$.
\endproclaim

\demo{Proof} Suppose that $a$ is $\alpha$-integrable in Exel's sense.  Exel
points
out (before $6.2$ of \cite{E1}) that $a$ is then integrable in the sense used in
Olesen and Pedersen discussed above (though they only consider positive
elements).  So if
$a \in A^+$ then we can apply Theorem 4.3 to conclude that $a \in {\Cal
P}_{\alpha}^+$.

Conversely, if $a \in {\Cal P}_{\alpha}^+$ then for every $b \in A$ the net
$\{bp_\lambda(a): \lambda \in {\Cal B}\}$ is norm Cauchy by definition, and
similarly for
$b$ on the other side.  The subnet obtained by restricting the $f$'s to be
characteristic functions of precompact subsets of $G$ is clearly cofinal,
so this
subnet too is Cauchy.  Similarly for $b$ on the other side.  Thus $a$ is
$\alpha$-integrable in Exel's sense. \qed
\enddemo

We remark that it follows by linearity that any element of ${\Cal
P}_{\alpha}$ is
$\alpha$-integrable in Exel's sense.

We tentatively make the following definition, which is the main one of this
paper. The reason that this definition is tentative will be explained in 
section 6.

\definition{4.5 Definition} The action $\alpha$ of $G$ on $A$ is {\it
proper} if
${\Cal P}_{\alpha}$ is dense in $A$.
\enddefinition

We will see later in Proposition 6.8 that
every action $\a$ has a canonical proper part, namely its restriction to the
closure of ${\Cal P}_\a$.

We show now that all of the 
examples successfully treated by the definition of
\cite{Rf7} are examples of proper actions in the sense of Definition
4.5. (See \cite{Rf8} and \cite{M} for further 
such examples in addition to those
already described
in \cite{Rf7}.)
This already gives a substantial supply of interesting examples.  The main
theorem
of \cite{E1} provides yet a further class of examples, associated to
$C^*$-algebraic
bundles over locally compact Abelian groups, to which Definition 4.5
applies.

\proclaim{4.6 Proposition} If the action $\alpha$ of $G$ on $A$ is proper
in the
sense of definition $1.2$ of \cite{Rf7}, then it is proper in the sense of
Definition 4.5 above.
\endproclaim

\demo{Proof} We recall that if $\alpha$ is proper in the sense of
definition $1.2$
of \cite{Rf7}, then there is a dense $*$-subalgebra $A_0$ of $A$ such that
if $a,b
\in A_0$ then the function $x \mapsto a\alpha_x(b)$ is integrable on $G$,
and for
$a,b \in A_0$ there is an $m \in M(A)^{\alpha}$ such that for every $c \in
A_0$ we
have
$$
\int c\alpha_x(a^*b)dx = cm.
$$
(There is more to the definition, but this suffices for our present purposes.)

We claim that $A_0 \subseteq Q_{\alpha}$, 
so that $A_0^2 \subseteq{\Cal P}_{\alpha}$.
Since $A_0^2$ (linear span of products) is dense in $A$ because $A_0$ is,
it will
follow that $\alpha$ is proper in the sense of Definition $4.5$.

So suppose that $a \in A_0$.  By hypothesis there is an $m \in
M(A)^{\alpha}$ such
that
$$
\int c\alpha_x(a^*a)dx = cm
$$
for every $c \in A_0$.  For a given $c \in A_0$ the function $x \mapsto
c\alpha_x(a^*a)$ is by assumption integrable on $G$, and so we can find an
increasing sequence $\{\lambda_n\}$ in ${\Cal B}(G)$ such that
$\{\lambda_n(x)c\alpha_x(a^*a)\}$
converges pointwise to $c\alpha_x(a^*a)$.  By the Lebesgue dominated
convergence
theorem, $c\int \lambda_n(x)\alpha_x(a^*a)dx$ converges in norm to $cm$.  Then
$cp_{\lambda_n}(a^*a)c^*$ increases up to $cmc^*$ in norm.  
Thus $m \ge 0$ since
$A_0$ is dense.  Furthermore $cp_{\lambda_n}(a^*a)c^* \le cmc^*$ for all $c \in
A_0$, and it follows that $p_{\lambda_n}(a^*a) \le m$.  Since our sequence
$\{\lambda_n\}$
can include any given element of ${\Cal B}$, it follows that
$p_\lambda(a^*a) \le m$
for all $\lambda \in {\Cal B}$, so that $a^*a$ is $\alpha$-integrable.

Finally, it is easily seen that the collection of $c$'s in $A$ for which
$cp_\lambda(a^*a)c^*$ converges to $cmc^*$ is norm closed.  But it contains
$A_0$,
and thus it is all of $A$.  Hence, the net $\{p_\lambda(a^*a)\}$ converges up
to $m$ in
the weak-strict topology.  But we saw in Proposition 
$2.4$ that this implies that
$p_\lambda(a^*a)$ converges strictly to $m$. \qed
\enddemo
 
We now want to show that if $A$ is commutative, then our definition of proper
action on $A$ captures the usual notion of proper action on a space.  This is a
somewhat subtle matter, as seen by examining Example $1.18$.  In fact, already
Green's original example \cite{G} will do --- he 
was concerned with closely related
matters.  His example is the case of Example $1.18$ in which $\{L_n\}$ is the
constant sequence $L_n = 1$.  In this case $\alpha$ is an integrable
action.  But
for $f$ as constructed in Example $1.18$ we have
$$
\int \alpha_t(f)dt = \cases
1 &\text{on ${\Cal O}_*$} \\
2 &\text{on ${\Cal O}_n$ for $n \ge 1$,}
\endcases
$$
which is not continuous on $M/\alpha$.  Thus $\alpha$ is not proper as an
action on
$C_{\infty}(M)$, much as Green \cite{G} showed that $\alpha$ is not proper
as an
action on $M$.  (As Green suggests there \cite{G}, the study of the
transformation
group $C^*$-algebras for actions of the kind described in Example $1.18$
might be
of some interest.)

\proclaim{4.7 Theorem} Let $\alpha$ be an action of a locally compact
group $G$ on
a locally compact space $M$, and so on $A = C_{\infty}(M)$.  Then
$\alpha$ as action on $A$ is proper in the sense of Definition $4.5$ if
and only if
$\alpha$
as an action on $M$ is proper.
\endproclaim

\demo{Proof} If $\alpha$ on $M$ is proper, then {}from the discussion at the
beginning of Section 2 it follows that $C_c(M)$ consists of $\alpha$-proper
elements, so that $\alpha$ on $C_{\infty}(M)$ is proper.  Equivalently,
definition
$1.2$ of \cite{Rf7} applies, so we can invoke Proposition $4.6$.

Suppose, conversely, that $\alpha$ is proper as action on $A$.  We show
that then
$\alpha$ is proper as action on $M$.  We can assume that the
$\alpha$-orbits in $M$
are closed, and that the stability subgroups are compact, for otherwise
$\alpha$ on
$A$ is not even integrable, 
by Proposition $1.17$.  Since ${\Cal P}_{\alpha}$ is
assumed dense, and is an ideal in this commutative case, it contains $C_c(M)$.

Let us show first that $M/\alpha$ is Hausdorff. As mentioned in section 1,
this does not follow {}from integrability of $\alpha$.
Let $m,n \in M$ and suppose they are in different $\alpha$-orbits.  Since
the orbit
$\alpha_G(n)$ is closed, we can find $f \in C_c(M)^+$ such that $f(m) > 0$
while
$f(\alpha_G(n)) = 0$.  Since $f \in {\Cal P}_{\alpha}$, 
$F = \int \alpha_x(f)dx$
exists and is continuous.  Clearly $F(m) > 0$ while $F(n) = 0$.  Thus $F$ is a
continuous function on $M/\alpha$ which separates $m$ and $n$.  Since $m$
and $n$
are arbitrary, it follows that $M/\alpha$ is Hausdorff.  So we assume {}from
now on
that $M/\alpha$ is Hausdorff.

Suppose now that $\alpha$ on $M$ is not proper.  It is easily seen {}from the
definition that there is then a compact subset, $K$, of $M$ such that $\{x
\in G:
\alpha_x(K) \cap K = \emptyset\}$ is not precompact.  Thus we can choose a
net
$\{k_\mu\}$ of elements of $K$ and a net $\{x_\mu\}$ of elements of $G$
such that
$\alpha_{x_\mu}(k_\mu) \in K$ for each $\mu$, but the net $\{x_\mu\}$ is not
precompact.  By the compactness of $K \times K$ we can find a subnet
$\{(x_{\nu},k_{\nu})\}$ of the net $\{(x_\mu,k_\mu)\}$ such that
$k_{\nu} \to k_0$ and $\alpha_{x_{\nu}}(k_{\nu}) \to k'_0$ for
points
$k_0$ and $k'_0$ of $K$.  Since $M/\alpha$ is Hausdorff, it follows
that
$k_0$ and $k'_0$ are in the same $\alpha$-orbit, so there is a $y \in G$
with $k'_0
= \alpha_y(k_0)$.  If we replace each $x_{\nu}$ by
$(y^{-1}x_{\nu})^{-1}$,
we find that $\alpha_{x_{\nu}^{-1}}(k_{\nu}) \to k_0$.

Choose $f \in C_c(M)^+$ such that  $\int f(\alpha_{y^{-1}}(k_0))dy = 1$.  
Since the orbit of $k_0$ is
closed, it meets the support of $f$ in a compact set.  Since the stability
subgroup of
$k_0$ is compact, we can find a compact subset $C$ of $G$ such that
$$
\int_C f(\alpha_{y^{-1}}(k_0))dy = \int_G f(\alpha_{y^{-1}}(k_0))dy=1.
$$
Let $\chi$ denote the characteristic function of $C$, so $\chi \in
{\Cal B}$.
Then
$$
(p_{\chi}(f))(k_0) = \int_C f(\alpha_{y^{-1}}(k_0))dy = 1.
$$
Now $p_{\chi}(f)$ is continuous, and so
$(p_{\chi}(f))(k_{\nu}) \to 1$ and
$(p_{\chi}(f))(\alpha_{x_{\nu}^{-1}}(k_\nu) \to 1$.  But
$$
(p_{\chi}(f))(\alpha_{x_{\nu}^{-1}}(k_\nu)) = \int_C
f(\alpha_{(x_{\nu}y)^{-1}}(k_{\nu}))dy 
= \int_{x_{\nu}C} f(\alpha_{y^{-1}}(k_{\nu}))dy.
$$
Since the net $\{x_{\nu}\}$ is not precompact, it is not eventually
contained
in $CC^{-1}$.  So we can find a subnet, which we still denote by
$\{x_{\nu}\}$,
such that $x_{\nu} \notin CC^{-1}$ for all $\nu$.  Then $C$ and
$x_{\nu}C$ are disjoint for each $\nu$, so
$$
\int_G f(\alpha_{y^{-1}}(k_{\nu})) \ge \int_C
f(\alpha_{y^{-1}}(k_{\nu}))dy 
+ \int_{x_{\nu}C} f(\alpha_{y^{-1}}(k_{\nu}))dy ,
$$
which converges to $1 + 1 = 2$.  Thus
$$
\liminf \int f(\alpha_{y^{-1}}(k_{\nu}))dy \ge 2.
$$
Since $\int f(\alpha_{y^{-1}}(k_0))dy = 1$, we see that $\int
\alpha_y(f)dy$ is not
continuous on $M$, so $\alpha$ as action on $A$ is not proper, a
contradiction. \qed
\enddemo

We conclude this section by showing that the C*-weights in the
present context have slightly better continuity properties than
we encountered earlier. We first need:

\proclaim{4.8 Proposition} Let $\alpha$ be an action of $G$ on $A$, and let
$\lambda \in
{\Cal B}$.  For any bounded approximate identity $\{e_{\nu}\}$ for $A$,
the net
$\{p_\lambda(e_{\nu})\}$ converges strictly in $M(A)$ 
to $(\int \lambda(x)dx)1$.
\endproclaim

\demo{Proof} 
We first remark that if
$h$ is a continuous function {}from a compact space
$K$ to
$A$, then for any $\varepsilon > 0$ there is a $\nu_0$ such that $\|h(x) -
e_{\nu}h(x)\| < \varepsilon$ for all $x \in K$ and all $\nu >
\nu_0$.
The same is true for $e_{\nu}$ on the right of $h(x)$.
This follows by using the compactness of $K$ to approximate
$h$ by a
finite sum $\sum \varphi_jh(x_j)$ where $\{\varphi_j\}$ is a suitable
partition of
unity on $K$. 
Now let $K$ denote the (compact) support of $f$.  Let
$\varepsilon > 0$ and $c \in A$ be given.  Then
$$
\align
cp_\lambda(e_{\nu}) - c(\int \lambda(x)dx) &= 
\int \lambda(x)(c\alpha_x(e_{\nu}) -
c)dx \\
&= \int \lambda(x)\alpha_x(\alpha_{x^{-1}}(c)e_{\nu} - \alpha_{x^{-1}}(c))dx.
\endalign
$$
When we apply the above remark to the function $x \mapsto \alpha_{x^{-1}}(c)$, 
we see
that we can find $\nu_0$ such that
$$
\|cp_\lambda(e_{\nu}) - c(\int \lambda(x)dx)\| < 
\varepsilon \text{ for } \nu >
\nu_0.
$$
Taking adjoints, we obtain the corresponding result for $c$ on the other
side. \qed
\enddemo

\definition{4.9 Definition} A completely positive map $p$ {}from $B$ to $A$
is said
to be {\it non-degenerate} if for some bounded approximate identity
$\{e_{\lambda}\}$ for $B$ the net $\{p(e_{\lambda})\}$ converges strictly
to $r1 \in M(A)$
for some $r \in {\Bbb R}^+$.
\enddefinition

This is just the definition at the top of page 49 of \cite{La} 
except that we do not require that $\|p\| = 1$. Because we here
require that $p$ be {\it completely} positive, we can apply some
of the results in \cite{La}. In particular, 
by Lemma $5.3$ of \cite{La} we will have $r = \|p\|$.

Notice now that Proposition $4.8$ states that each $p_\lambda$ is 
non-degenerate, for $\lambda \in {\Cal B}$.

\definition{4.10 Definition} We will say that a strict $C^*$-weight $\psi$
{}from $B$
toward $A$ is {\it non-degenerate} if there is an increasing net
$P=\{p_{\lambda}\}$
of completely positive maps {}from $B$ to $A$ as in Definition $2.8$ such that
eventually each $p_{\lambda}$ is non-degenerate.
\enddefinition

According to Corollary $5.7$ of \cite{La}, if $p$ is non-degenerate, then $p$
extends uniquely to a completely positive map ${\bar p}$ {}from $M(B)$ to
$M(A)$ such
that ${\bar p}(1_{M(B)}) = \|p\|1_{M(A)}$, and ${\bar p}$ is strictly
continuous on
bounded subsets of $M(B)$.  This makes possible the following strengthening
of the
lower semi-continuity property stated in Proposition $2.9$, when $\psi$ is
non-degenerate.

\proclaim{4.11 Proposition} (Compare with $3.5$ of \cite{Ku}.)  Let $\psi$ be a
strict $C^*$-weight {}from $B$ toward $A$, and assume that $\psi$ is
(completely positive and) non-degenerate.
Let $b \in {\Cal P}_P^+$, and let $\{b_{\mu}\}$ be a 
net in $B^+$ which converges
strictly to $b$ and is such that $b_{\mu} \le b$ for each $\mu$.  Then the net
$\{\psi(b_{\mu})\}$ converges strictly to $\psi(b)$.
\endproclaim

\demo{Proof} The proof is the same as that for Proposition $2.9$ except
that now we
use the strict continuity of $p_{\lambda}$ in order 
to choose $\mu_0$ such that for $\mu \ge \mu_0$
$$
\|ap_{\lambda}(b)a^* - ap_{\lambda}(b_{\mu})a^*\| < \varepsilon/2.
$$
\qed
\enddemo

\section{5. Functoriality, and C*-algebras proper over a space}

In considering functoriality it is useful for us 
to treat ``morphisms'' \cite{La, Wr}. Let $A$ and $B$ 
be C*-algebras. A morphism {}from $B$ to $A$ is a homomorphism $\theta$
{}from $B$ to $M(A)$ which is non-degenerate in the sense that 
$\theta(B)A$ is dense in $A$. Then $\theta$ extends uniquely to a homomorphism,
$\bar \theta$, {}from $M(B)$ to $M(A)$, which is strictly continuous on
bounded sets \cite{La}. If $\a$ and $\beta$ are actions of $G$ on $A$ and
$B$, then we say that $\theta$ is equivariant if 
$$
   \theta(\beta_x(b)) = \a_x(\theta(b))
$$
for all $x \in G$ and $b \in B$, where here $\a$ has been extended
to $M(A)$. The extension of $\theta$ to $M(B)$ is
then seen to be equivariant in the usual sense. The following 
proposition is basically proposition 1.4 of \cite{QR} once Theorem 2.6
above is taken into account.

\proclaim{5.1 Proposition} With $\a$, $\beta$ and $\theta$ as above, we
have $\bar \theta({\tilde {\Cal P}}_\beta) \subseteq {\tilde {\Cal P}}_\alpha$,
and 
$$
     {\tilde \psi}_\a(\bar \theta(n)) = \bar \theta({\tilde \psi}_\beta(n))
$$
for all $n \in {\tilde {\Cal P}}_\beta$.
\endproclaim

\demo{Proof} It is easily seen that $\bar \theta(M(B)_e) \subseteq M(A)_e$.
Let $n \in {\tilde {\Cal P}}_\beta$. By definition the 
net $\{p_\lambda^\beta(n)\}$ is bounded and converges 
to ${\tilde \psi}_\beta(n)$ strictly. 
Thus $\{{\bar \theta}(p_\lambda^\beta (n))\}$ converges strictly 
to ${\bar \theta}({\tilde \psi}_\beta(n))$. 
But ${\bar \theta}(p_\lambda^\beta(n)) = p_\lambda ^\a({\bar \theta}(n))$.
\qed
\enddemo

The next lemma should be compared carefully with 
the definition of hereditary (non-closed) subalgebras
in VII.4.1 of \cite{FD}.

\proclaim{5.2 Lemma} Let $\Cal H$ be any hereditary $\ast$-subalgebra, 
possibly not
closed, of a C*-algebra $C$. Then ${\Cal H}C{\Cal H} \subseteq {\Cal H}$.
Let $D$ be the closure 
$({\Cal H}C{\Cal H})\ \bar {}$, where  ${\Cal H}C{\Cal H}$
means linear span. Then  ${\Cal H}C{\Cal H} \cap C^+$ is dense
in $D^+$.
Furthermore, the closure, $\bar{\Cal H}$, of ${\Cal H}$ in $C$ is a hereditary 
C*-subalgebra of $C$.
\endproclaim

\demo{Proof} Here, in contrast to \cite{FD}, by 
``hereditary'' we mean that $\Cal H$ is the linear 
span of its positive part ${\Cal H}^+$, and that ${\Cal H}^+$ is a hereditary
cone in $C$ in the sense we used earlier. 
Let $c \in C^+$ and $h \in \Cal H$. Then
$h^*ch \leq \|c||h^*h$, so that $h^*ch \in \Cal H$. By linearity
it follows that this is true for any $c \in C$. By polarization it then
follows that if $h_1, h_2 \in \Cal H$ and $c \in C$, then $h_1ch_2
\in \Cal H$.  

Now suppose that $d \in D^+$. By considering an approximate identity
for $D$, and approximating its elements by elements of  ${\Cal H}C{\Cal H}$,
we can approximate $d$ by elements of form $b^*db$ where
$b \in  {\Cal H}C{\Cal H}$. But 
then $b^*db \in ({\Cal H}C{\Cal H}) \cap C^+$. 

Finally, it is clear that $\bar{\Cal H}C\bar{\Cal H} \subseteq \bar{\Cal H}$.
But as indicated in VII.4.1 of \cite{FD} this does imply that $\bar{\Cal H}$
is hereditary in our sense, since $\bar{\Cal H}$ is closed.     \qed
\enddemo

It is easily seen that if there is an equivariant map {}from a $G$-space $Y$
to a $G$-space $Z$ and if the action on $Z$ is proper, then the action
on $Y$ must be proper. We have the following generalization to the 
non-commutative case:

\proclaim{5.3 Theorem} Let $\a$ and $\beta$ be actions of $G$ on C*-algebras
$A$ and $B$, and suppose that $\beta$ is proper. If there is an
equivariant morphism {}from $B$ to $A$, then $\a$ is proper.
\endproclaim

\demo{Proof} Let $\theta$ be an equivariant morphism {}from $B$ to $A$.
Since $\theta$ is non-degenerate, $\theta(B)A\theta(B)$ is dense
in $A$. Since $\beta$ is proper, ${\Cal P}_\beta$ is dense in $B$, and
so $\theta({\Cal P}_\beta)A\theta({\Cal P}_\beta)$ is dense in $A$. 
But $\theta({\Cal P}_\beta) \subseteq {\tilde {\Cal P}}_\a$ by Proposition 5.1.
Thus ${\tilde {\Cal P}}_\a A {\tilde {\Cal P}}_\a$ is dense in $A$.
But ${\tilde {\Cal P}}_\a$ is a hereditary $\ast$-subalgebra in $M(A)$ by
Key Lemma 2.7. Thus ${\tilde {\Cal P}}_\a A {\tilde {\Cal P}}_\a 
\subseteq {\tilde {\Cal P}}_\a$
by Lemma 5.2. Since ${\Cal P}_a = {\tilde {\Cal P}}_\a \cap A$,
it follows that ${\Cal P}_\a$ is dense in $A$. 
    \qed
\enddemo

\proclaim{5.4 Corollary} Let $\a$ be a proper action of $G$ on
a C*-algebra $A$, and let $I$ be an $\a$-invariant ideal in $A$, so
that $\a$ drops to an action, $\bar \a$, on $A/I$. Then $\bar \a$
is proper.
\endproclaim

\proclaim{5.5 Proposition} For $\a$, $A$ and $I$ as just above, the action 
defined by $\a$ on $I$ is proper.
\endproclaim

\demo{Proof} Since ${\Cal P}_\a$ is dense in $A$, the linear span
${\Cal P}_\a I {\Cal P}_\a$ must be dense in $I$. But it is contained 
in ${\Cal P}_\a$ by Lemma 5.2. Thus ${\Cal P}_\a \cap I$ is dense in $I$. In
fact, {}from Lemma 5.2 we see that ${\Cal P}_\a \cap I^+$ is
dense in $I^+$. Each element of $M(A)$ determines an element of $M(I)$
in the evident way. Let $c \in {\Cal P}_\a \cap I^+$, and let $\psi_\a(c)$
also denote the corresponding element of $M(I)$. It is easily seen that
$\{p_\lambda(c)\}$ converges strictly to $\psi_\a(c)$ in $M(I)$. \qed
\enddemo

An increasingly important way of getting aspects of properness
to bear on an action of a group on a C*-algebra is to have an
equivariant morphism {}from a commutative C*-algebra with proper action,
whose image is central. This technique seems to have been first introduced
by Kasparov, in section 3 of \cite{Ks}. For more recent occurences
see \cite{GHT} and the references therein. Such a morphism is a
special case of the situation of Theorem 5.3, so that we
immediately obtain:

\proclaim{5.6 Corollary} Let $\a$ be an action of $G$ on a
C*-algebra $A$. Let $\beta$ be a proper action of $G$ on a locally
compact space $Z$, and so on $C_\infty(Z)$. If there is an equivariant
morphism {}from $C_\infty(Z)$ to $A$ whose image is contained in the
center of $M(A)$, then $\a$ is proper.
\endproclaim

The deficiency of this corollary is that, as we discuss in the next
section, we have not seen how to establish strong Morita equivalence
between the generalized fixed-point algebra and an ideal in the
crossed product in the general situation of our present definition
of properness. However, we now show that, even in the absence of
the requirement that the image of the morphism be central, 
the situation of Corollary
5.6 falls within the perview of definition 1.2 of \cite{Rf7}, where
we were able to establish this Morita equivalence. The fact that
centrality is not needed seems to be a new observation. Our proof can be
viewed as a variation of the proof of theorem 3.13 of \cite{Ks},
with some ingredients also {}from \cite{Qg}. We will
not include discussion of the fact that if the action on $Z$ is also
free, then the strong Morita equivalence is with the whole crossed
product, but see the discussion of ``saturation'' in \cite{Rf7}.

\proclaim{5.7 Theorem} Let $\a$ be an action of $G$ on a C*-algebra $A$.
Let $\beta$ be a proper action of $G$ on a locally
compact space $Z$, and so on $C_\infty(Z)$, and let $\theta$ be an
equivariant morphism {}from $C_\infty(Z)$ to $A$. 
Let $A_0$ denote the linear span of
$(\theta(C_c(Z))A(\theta(C_c(Z))$, which is a dense $\ast$-subalgebra of $A$. 
Then $A_0$ satisfies the conditions of definition 1.2 of \cite{Rf7}, so
that $\a$ is proper in the sense of that definition. Thus the generalized
fixed-point algebra (in the sense of \cite{Rf7}) is strongly Morita
equivalent to an ideal in the reduced crossed product algebra.
\endproclaim

\demo{Proof} For notational simplicity we sometimes omit $\theta$ 
below, and confuse $\beta$ with $\a$. Let $a, b \in A$ and
$f, g \in C_c(Z)$, and consider the function
$$
  x \mapsto (af)\a_x(gb) = a(f\beta_x(g))\a_x(b).
$$
Since $\beta$ is proper, this function has compact support. 
{}From this it is easily seen that if $a, b \in A_0$, then the
function $x \mapsto a\a_x(b)$ has compact support, and
so is in $L^1(G, A)$, as will be its product with $\Delta ^{-1/2}$.
This says exactly that condition 1 of definition 1.2 of \cite{Rf7} is
satisfied.

We turn now to condition 2. By essentially the same argument as in
the proof of Theorem 5.3, using the fact 
that $C_c(Z) \subseteq {\Cal P}_\beta$, we find 
that $A_0 \subseteq {\Cal P}_\alpha$.
For the element $\<a, b\>_D$ of $M(A)^\a$ which is required 
by condition 2 for any $a, b \in A_0$ we take $\psi_\a(a^*b)$. (See
Proposition 1.12 for the $\a$-invariance.)
Then for $c \in A_0$ we have
$$
  c\<a, b\>_D = \lim \int \lambda(x)c\a_x(a^*b)dx.
$$
But as seen above, $x \mapsto c\a_x(a^*b)$ has compact support,
and so the net of integrals is eventually constant, and has limit
$$
 \int c\a_x(a^*b)dx,
$$
as required by condition 2. Now condition 2 also requires that
$c\<a, b\>_D$ be again in $A_0$ for $a, b, c \in A_0$, that is,
that $\<a, b\>_D \in M(A_0)^\a$ in the notation of \cite{Rf7}. 
It is easily seen that it suffices to show that if $a, b \in A$
and if $f_1, f_2, g_1, g_2 \in C_c(Z)$, then
$f_1af_2\psi_\a(g_1bg_2) \in A_0$. But by the argument {}from
the proof of condition 1 above we see that this element is
given by
$$
  f_1\int a(f_2\a_x(g_1))\a_x(b)\a_x(g_2)dx  .
$$
Let $K$ denote the support of the integrand, which is compact.
Let $S$ denote the support of $g_2$, and let $L = \a_K(S)$. Choose
$h \in C_c(Z)$ such that $h \equiv 1$ on $L$. Then for $x \in K$
we have $\a_x(g_2) = \a_x(g_2)h$. Consequently the above 
expression can be rewritten as
$$
 f_1(\int a(f_2\a_x(g_1))\a_x(b)\a_x(g_2)dx)h,
$$
which is manifestly in $A_0$. 
\qed
\enddemo

It is not clear to me how Thomsen's definition of a K-proper action, given
in 9.1 of \cite{T}, relates to our present considerations, though it
has some relation to \cite{GHT}.

\section{6. Strong Morita equivalence}

In this section we will discuss what one might take as
the ``generalized fixed-point algebra" for a proper
action. Our guiding principle will be that this
generalized fixed-point algebra should be strongly Morita
equivalent \cite{Rf2, Rf4, Rf5} to at least an ideal in the crossed product
algebra, much as happens in \cite{Rf7}. (In the case of
commutative $A=C_{\infty}(X)$  we know \cite{Rf7} that
it will be strongly Morita
equivalent to the entire reduced crossed product
algebra exactly if the action on $X$ is free.) The outcome of our
discussion will be far {}from satisfactory. In particular,
Exel \cite{E2} has shown that the candidate for ``generalized fixed-point algebra''
which I had suggested in the first version of this paper is often
too big. (See our discussion following Theorem 8.5.) In fact, Exel \cite{E2} has
shown that the situation is fairly subtle, and the question of how best to define
the generalized fixed-point algebra remains unresolved.

We let $M(A)^{\alpha}$ denote the subalgebra of
elements in $M(A)$ which are fixed by $\alpha$.  {}From Proposition $1.12$ we
immediately obtain:

\proclaim{6.1 Proposition} The range of $\psi_{\alpha}$ is contained in
$M(A)^{\alpha}$.
\endproclaim

By viewing elements of $M(A)^\alpha$ as constant functions on $G$ and
applying Proposition 3.1 we obtain:

\proclaim{6.2 Proposition} Let $a \in {\Cal P}_{\alpha}$ and let $m \in
M(A)^{\alpha}$.  Then $ma$ and $am$ are in ${\Cal P}_{\alpha}$, and
$$
\psi_{\alpha}(ma) = m\psi_{\alpha}(a),\hskip 0.5 in \psi_{\alpha}(am) =
\psi_{\alpha}(a)m.
$$
If $a \in Q_{\alpha}$ then $am \in Q_{\alpha}$.
\endproclaim

We note that this proposition implies that the range of $\psi$ is an ideal
in $M(A)^{\alpha}$ (and clearly a $*$-ideal).

It is clear {}from Proposition $3.1$ that ${\Cal P}_{\alpha}$ and ${\Cal
Q}_{\alpha}$ are right $M(A)^\a$-modules.  Then $\psi_{\alpha}$ is almost a
generalized conditional expectation {}from ${\Cal P}_{\alpha}$, as defined in
definition $4.12$ of \cite{Rf2}.  The only property which is not clear is the
$\psi$-density of ${\Cal P}_{\alpha}^2$ in ${\Cal P}_{\alpha}$ as defined in
property 5 of definition $4.12$ of \cite{Rf2}.  I do not know how often it
holds.
(The relative boundedness of property 4 of the definition follows {}from the fact
that for $b \in {\Cal P}_{\alpha}$ the map $a \mapsto \psi_{\alpha}(b^*ab)$ is
defined on all of $A$ and is positive, and so is bounded.)

We remark that the KSGNS construction of \cite{Ku} can, of course, be
carried out
in our special case.  Because of the above conditional expectation property of
$\psi_{\alpha}$, the KSGNS construction for $\psi_{\alpha}$ is in
this case essentially the
``induction in stages'' construction of theorem $6.9$ of \cite{Rf2}, applied to
${\Cal Q}_{\alpha}$ as left-$A$ right-Hilbert-$M(A)^\a$-module using
$\psi_{\alpha}$, and $A$ as left-$M(A)^\a$ right-Hilbert-$A$-module in the
canonical way.  By $3.7$ of \cite{Ku} the construction gives a non-degenerate
representation of $A$.  The proof of non-degeneracy basically uses Proposition
$2.12$ above.

It is at first sight not entirely clear what one should
take as the ``generalized fixed-point algebra". Our
guiding principle will be our desire, just expressed above,
that it be  strongly Morita
equivalent to at
least an ideal in the crossed product algebra. One
possibility is to take the generalized fixed-point
algebra to be the closure of the range of
$\psi_{\alpha}$. We now give an example to show that
already when $A$ is commutative this does not accord with
our guiding principle.

\example{6.3 Example} Let $G$, $A$ and $\alpha$ be as in
Example 2.10. Then it is easily seen that
$A\times_{\alpha}G$ is isomorphic to the
$C^*$-direct-sum of a countable number of compact
operator algebras. Its primitive ideal space is thus a
countable discrete set, and so it cannot be strongly
Morita equivalent to a unital $C^*$-algebra, since
strongly Morita $C^*$-algebras have homeomorphic
primitive ideal spaces (corollary 3.3 of \cite{Rf3}). But
let $g$ be as in Example 2.10. It is seen there that
$\psi(g)=1$. So the closure of the range of $\psi$ is a
unital $C^*$-algebra, and thus cannot be strongly
Morita equivalent to $A\times_{\alpha}G$. Of course the
difficulty is that in this case we want the generalized
fixed-point algebra to be contained in $C_\infty(N)$.
\endexample

To try to remedy the situation we consider a slightly subtler definition.
Since $P_\a=Q^*_\a Q_\a$, we can define an
$M(A)^\a$-valued inner-product $\< \ , \ \>_D$, on
$Q_\a$ by
$$
\< a,b\>_D=\psi(a^*b) \ .
$$
By Proposition 3.1 this behaves correctly for the right
action of $M(A)^\a$. But since $P_\a=Q^*_\a Q_\a$, the
span of the range of this inner-product is just the
range of $\psi_{\a}$, and so by the above example this span
will not be appropriate as the generalized fixed-point
algebra.  So instead, we consider the restriction of this
inner-product to $P_\a\subseteq Q_\a$. Exel's example \cite{E2}
shows that this is in general still too big. (See the discussion
after Theorem 8.5 below.) But we
examine it here. That is, we set:

\definition{6.4 Definition} Let $\a$ be an action of
$G$ on $A$. The {\it big generalized fixed-point} algebra of $\a$ is
the norm closure of the linear span of the elements of $M(A)^{\a}$ of 
the form
$$
\< a,b\>_D=\psi_{\a}(a^*b) \ 
$$
for  $a,b\in P_{\a}$. We will denote it by $D_{\a}$.
\enddefinition

This accords with the definition given in \cite{Rf7},
as well as with definition 1.5 of \cite{QR}. It is clear {}from
Proposition 6.2 that $D_\alpha$ is an ideal in $M(A)^\alpha$.
We remark that just for the purpose of stating this definition we do
not need to assume that $\a$ is proper. 

We now proceed to show that at least when $A$ is commutative
this definition provides a generalized fixed-point algebra
where we want it.

\proclaim{6.5 Proposition} Let $\alpha$ be a proper action of $G$
on a commutative C*-algebra $A = C_\infty(M)$. Let $M/\alpha$ be the
orbit space, with its quotient locally compact Hausdorff topology.
For $f, g \in {\Cal P}_\alpha$ we have
   $$\psi_\alpha(\bar{f}g) \in C_\infty(M/\alpha).$$
\endproclaim

\demo{Proof}
Of course $\psi_\alpha(\bar{f}g) \in C_b(X/\alpha)$. The only issue
is the vanishing at infinity. Let $\e \geq 0$ be given. We can
find compact $K \subseteq M$ such that $|f(m)| \leq \e$ for $m \notin K$. 
The image, $\dot K$, of $K$ in $M/\alpha$ is compact. Let $m \in M$
with $\dot m \notin \dot K$. View evaluation at $m$ as a state of $A$, so
that we can apply Theorem 3.4. Then $x \mapsto
(\bar f g)(\alpha_x^{-1}(m))$ is integrable on $G$, and
  $$|(\psi(\bar f g))(m)| \leq 
       \int_G|\bar f (\alpha_x^{-1}(m)g(\alpha_x^{-1}(m))| dx
      \leq \e\int_G|g(\alpha_x^{-1}(m)| dx \leq\e\|\psi_\alpha(g)\|_\infty.
$$
That is, $|\<f,\ g\>_D(\dot m)| \leq \e\|\psi_\alpha(g)\|_\infty$ 
for $\dot m \notin \dot K$, as
desired.       \qed
\enddemo

\bigskip We now turn to the question of strong Morita equivalence. It
is natural in view of \cite{Rf7} to take $\Cal P_{\a}$ (suitably
completed) as the equivalence bimodule. We know that it is a right
$D_{\a}$-module. We restrict to ${\Cal P}_\alpha$ the inner-product 
defined above on ${\Cal Q}_\alpha$. It then has values in
$D_{\a}$ by definition. We consider the
corresponding norm $\|a\|_{\a}=\|\<a,a\>_D\|^{\frac 12}$.
Then the completion, $\bar\Cal P_{\a}$, for this norm is a right
Hilbert $D_{\a}$-module \cite{La}, which is {\it full} in the sense
that the span of the range of the inner-product is dense in $D_{\a}$
(because we defined $D_{\a}$ that way).

Thus we have the corresponding algebra, $B(\bar\Cal P_{\a})$, of
bounded (adjointable) operators on $\bar\Cal P_{\a}$, and its ideal
$E$ of ``compact" operators on $\bar\Cal P_{\a}$ \cite{Rf2, La}
generated by the ``rank-one" operators $\<a,b\>_E$ defined by
$$
\<a,b\>_Ec = a\<b,c\>_D \ .
$$
Always $E$ is strongly Morita equivalent to $D_{\a}$ \cite{Rf2}.
What we need to do is to relate $E$ to the (reduced) crossed product
algebra $A\x^r_{\a}G$. So we examine the extent to which
$A\x^r_{\a}G$  acts on $\Cal P_{\a}$. We begin by considering the
action of $G$.

For each $x\in G$ we define an operator, $U_x$, on $\Cal P_{\a}$ by
the same formula as in Notation 1.15. That this operator carries
$\Cal P_{\a}$ into itself follows {}from the following more general
fact:

\proclaim{6.6 Proposition} Let $\mu$ be a finite measure of compact
support on $G$. For any $a\in\Cal P_{\a}$ define $U_{\mu}a$ by
$$
U_{\mu}a=\int_G U_xa \ d\mu(x) \ ,
$$
in terms of the $C^*$-norm of $A$. Then  $U_{\mu}a\in\Cal P_a$, and
$$
\psi_{\a}(U_{\mu}a) = (\int\Delta(y)^{-\frac 12} d\mu(y))\psi(a) \ .
$$
\endproclaim

\demo{Proof}  It suffices to prove this for the case in which $\mu$
is positive and $a\in\Cal P_{\a}^+$. Now for $\lambda\in\Cal B$ we have
$$
\align
p_{\lambda}(U_{\mu}a) &= \int \lambda(x)\a_x
 (\int \Delta (y)^{\frac 12}\a_y(a) \ d\mu(y))dx \\
&= \int \Delta (y)^{\frac 12} \int \lambda(x)\a_{xy}(a)dx \ d\mu(y)
= \int \Delta (y)^{\frac 12} \int \lambda(xy^{-1}) \a_x(a)
 \Delta(y^{-1})dx \ d\mu(y)\\
&=\int (\int  \lambda(xy^{-1}) \Delta (y)^{-\frac 12}  \ d\mu(y))
\a_x(a)dx \ .
\endalign
$$
Denote the inner integral by $\lambda *\mu$. It is in $C^+_c(G)$,
and we can rewrite the above as $p_\lambda(U_{\mu}a)=p_{(\lambda*\mu)}(a)$.
If we scale $\mu$ so that $ \Delta (y)^{-\frac 12}   d\mu(y)$ is a
probability measure, then $\lambda *\mu\in\Cal B$. Furthermore, the
collection of such $\lambda *\mu$'s is cofinal in $\Cal B$, since
given $\lambda_1\in {\Cal B}$ we can choose $\lambda$ such that it has value
1 on (support$(\lambda_1)$)(support$(\mu))^{-1}$, so that
$\lambda *\mu$ has value 1 on support$(\lambda_1)$. Consequently
$p_{\lambda}(U_{\mu}a)$ must converge strictly to $\psi(a)$.
Then in view of how we scaled $\mu$, we obtain the desired
conclusion. \qed
\enddemo

Note that if $\mu$ does not have compact support,
$\int\Delta(y)^{-\frac 12}d\mu(y)$ may not be finite.

\proclaim{6.7 \ Corollary} The action $\a$ carries $\Cal P_{\a}$ into itself.
\endproclaim

Let me remark
that I do not know whether $U_{\mu}$ carries $\Cal Q_{\a}$ into
itself in general.

We can now clarify a remark made after Definition 4.5.

\proclaim{6.8 Proposition} Let ${\bar {\Cal P}}_\a$ denote the
norm closure of ${\Cal P}_\alpha$ in $A$. Then ${\bar {\Cal P}}_\a$
is a hereditary C*-subalgebra of $A$ which is carried into itself
by $\a$, and on which the action $\a$ is proper.
\endproclaim

\demo{Proof} Denote ${\bar {\Cal P}}_\alpha$ by $B$. It is clear that
$B$ is a C*-subalgebra of $A$, which {}from Corollary 6.7 is carried
into itself by $\a$. Since ${\Cal P}_\a$ is hereditary, it follows
{}from Lemma 5.2 that $B$ is hereditary.

We now show that the action $\a$ on $B$ is proper. For clarity of 
argument we denote this restricted action by $\beta$. Let 
$a \i {\Cal P}_\a^+$. It suffices to show that then $a \in {\Cal P}_\beta^+$. 
Since $a \in {\Cal P}_\a^+$, there is an $m \in M(A)^+$ such that
the net $\{p_\lambda(a)\}$ converges strictly to $m$. It is easily
seen that $p_\lambda(a) \in B$ for each $\lambda$, since $a \in B$.
For $c \in B$ the net $\{cp_\lambda(a)\}$ is in $B$ and converges
in norm to $cm$. Thus $cm \in B$. Similarly $mc \in B$. That is,
$m$ normalizes $B$, and so determines an element, say $n$, of
$M(B)$. The above steps then show that $\{p_\lambda(a)\}$ converges strictly
for $B$ to $n$. Thus $a \in {\Cal P}_\beta^+$ as desired.     \qed
\enddemo

We remark that the above proposition does not adequately capture
the notion of the wandering subset of an action on an ordinary
space. For example, let $M$ be the one-point compactification of
$\Bbb R$, with action $\a$ of $G = {\Bbb R}$ by translation, leaving
the point at infinity fixed. The wandering subset is $\Bbb R$, on
which the action is proper. But if we set $A = C(M)$ with corresponding 
action $\a$, then it is easily checked that ${\Cal P}_\a = \{0\}$. 

Since $\psi_{\a}$ is the restriction to $\Cal P_{\a}$ of the
$\phi_{\a}$ of section 1, it follows as in Notation 1.15 that $U_x$
is ``unitary" for the $D_{\a}$-valued inner product on  $\bar\Cal
P_{\a}$.

However, we also need the representation $x\mapsto U_x$ of $G$ on $\bar\Cal
P_{\a}$ to be strongly continuous for the norm $\| \ \|_{\a}$ on $\bar\Cal
P_{\a}$. I have not succeeded in showing that this holds in general,
though it can be shown to hold in many examples. This kind of
question is known to be difficult even in the case of ordinary
weights (in contrast to traces). See lemma 3.1 of \cite{QV} for a
fairly restrictive hypotheses, ``regularity" (also discussed in
\cite{Ku}), under which one can
prove this strong continuity for weights.

It is natural to expect that $U$ is strongly continuous on vectors of
the form $U_ga$ where $a\in\Cal P_{\a}$ and $g\in C_c(G)$
and where we view $g$ as the measure $\mu =g(x)dx$. We now show that
this is the case. But this then reduces our question to:

\example{6.9 Question.}  With notation as above, is the linear
span of the elements of $\Cal P_{\a}$ of the form $U_ga$ for $a\in
\Cal P_{\a}$ and $g\in C_c(G)$, dense in $\Cal P_{\a}$ for the norm
{}from the $D_{\a}$-valued inner product?
\endexample

\proclaim{6.10 \ Lemma} Let $\mu$ be any finite complex measure of
compact support on $G$ and let $a,b\in\Cal P_{\a}$.  Then
$$
\|\<b,U_\mu a\>_D\| \leq 4\|\psi_\a(b^*b)\| \
\|\psi_\a(a^*a)\| \ \|\mu\|_1 \ ,
$$
where $\|\mu\|_1$ denotes the total variation norm of $\mu$.
\endproclaim    

\demo{Proof} Let $\lambda\in\Cal B$.  Then
by Proposition 3.3 and the ``unitarity" of $U$
$$
\|p_{\lambda}(b^*U_xa)\|  \leq 4\|\psi_\a(b^*b)\| \
\|\psi_\a(a^*a)\|  \ .
$$
Consequently
$$
\|p_{\lambda}(b^*U_{\mu}a)\|  \leq 4\|\psi_\a(b^*b)\| \
\|\psi_\a(a^*a)\|  \ \|\mu\|_1 \ .
$$
But $p_{\lambda}(b^*U_{\mu}a)$ converges strictly to
$\<b,U_{\mu}a\>_D$. The desired inequality follows.  \qed
\enddemo

\proclaim{6.11 Proposition} Suppose that $a\in\Cal P_{\a}$ is of
the form $U_g(d)$ for $d\in\Cal P_{\a}$ and $g\in C_c(G)$. Then the
function $x\mapsto U_xa$ is continuous for the norm on $\Cal P_{\a}$
coming {}from the $D_\a$-valued inner product defined by $\psi_\a$.
\endproclaim

\demo{Proof} Since $U_g(d)$ is defined in terms of the norm of $A$,
a standard calculation shows that $U_x(U_gd)=U_{L_xg}(d)$ where $L_x$
is the usual left translation on function on $G$. Consequently for
any $b\in\Cal P_{\a}$,
$$
\align
\|\<b,U_xa\>_D-\<b,a\>_D\| & =
\|\<b,U_xa-a\>_D\| \\
 = \|\<b, U_{(L_xg-g)}d\>_D\| & \leq
4\|\psi(b^*b)\|  \ \|\psi(d^*d)\|\| L_xg-g\|_1 \ .
\endalign
$$
But it is a standard fact that $L$ is strongly continuous on
$L^1(G)$. {}From this the above inequality shows that $U$ is ``weakly
continuous". The ``strong continuity" then follows in the usual way
{}from the fact that $U$ is ``unitary". That is,
$$
\align
\|U_a-a\|^2_{\psi} & =
\|\< U_xa-a, \ U_xa-a\>_D\| \\
&= \|\<a,a\>_D-\<U_xa,a\>_D-\<a,U_xa\>_D+\<a,a\>_D\| \\
&= 2\|\<a-U_xa, \ a\>_D\|
\endalign
$$
\qed
\enddemo

\bigskip
As our equivalence bimodule we should surely take the part of $\Cal
P_{\a}$ on which $U$ is strongly continuous for the norm {}from $\psi_\a$,
which the above proposition makes clear is still dense in $A$ if
$\a$ is proper. But since I have not seen how to overcome the main
obstacle, which we will discuss shortly, I will avoid the added
notational complexity this would require in view of the lack of an
answer to Question 6.8. We will just continue to deal with $\Cal
P_{\a}$ itself. Note also that if $G$ is discrete the issue of strong
continuity does not arise.

We now turn to the action of $A$. The left action, $L$, of
$P_{\a}$ on itself commutes with the right action of $D_\a$.
Let $b\in Q_{\a}$. Since $Q_{\a}$ is a left ideal in $A$, the
positive linear functional $a\mapsto \psi_{\a}(b^*ab)$ is defined on
all of $A$, and so is continuous (lemma 6.1 of \cite{La}). Thus
there is a constant, $K$, such that
$$
\|\psi_{\a}(b^*a^*ab)\|\leq K\|a\|^2 \ ,
$$
and this remains true when $a$ and $b$ are restricted to be in
$\Cal P_{\a}$. Then this says that
$$
\|\<L_ab, \ L_ab\>_D\|\leq K\|a\|^2 \ ,
$$
so that the $*$-homomorphism $L$ of
 $\Cal P_{\a}$ into $B(\bar\Cal P_{\a})$ is continuous, hence of norm
1. We have thus obtained:

\proclaim{6.12 Proposition}
For  $a\in\Cal P_{\a}$ the operator $L_a$ on $\Cal P_{\a}$ is a
bounded operator for the $D_\alpha$-valued inner product, and in fact
$\|L_a\|\leq \|a\|$. Thus $L_a$ extends to a continuous
$*$-homomorphism {}from the closure in $A$ of $\Cal P_{\a}$  into
$B(\bar\Cal P_{\a})$.
\endproclaim

However, in general I do not see why the representation $L$ on
$\bar\Cal P_{\a}$ need be non-degenerate, i.e.~why
$L_{\Cal P_{\a}}(\bar\Cal P_{\a})$ need be dense in $\bar\Cal P_{\a}$
for the norm {}from
$\psi_\a$, although again this can be shown to be true for many examples.
This question is closely related to the question mentioned
in the comments following Proposition 6.2.

On the other hand, $U$ and $L$ do satisfy the usual covariance
relation. For $a,b\in\Cal P_{\a}$ and $x\in G$ we have
$$
U_x(L_ab)=\Delta(x)^{\frac 12}\a_x(ab)=
L_{\a_x(a)}U_x \ b \ ,
$$
so that
$$
U_xL_a=L_{\a_x(a)}U_x \ .
$$
Thus if $U$ is strongly continuous (for example if $G$ is discrete)
and if the representation $L$ of $A$ is non-degenerate, then by the
usual universal property \cite{Pe2} we obtain a non-degenerate
$*$-representation of the crossed product algebra $A\times_{\a}G$ on
$\bar\Cal P_{\a}$. We will not repeat here the discussion {}from
\cite{Rf6} which indicates that we should actually obtain a
representation of the reduced crossed product algebra, since we have
more serious difficulties. Namely, we need to know that the algebra $E$ of
compact operators, which is strongly Morita equivalent to $D_{\a}$,
is contained in (the image of) the crossed product algebra. For this
it suffices to show that whenever $a,b\in\Cal P_{\a}$ then
$\<a,b\>_E$, defined above, is in $A\times_{\a}G$. Now at least
symbolically, for $c\in\Cal P_{\a}$,
$$
\<a,b\>_Ec=a\<b,c\>_D=\int a\a_x(b^*)\a_x(c)dx
= \int a\a_x(b^*)\Delta(x)^{-\frac 12} U_x(c)dx \ .
$$
So we want the function $x\mapsto a\a_x(b)\Delta(x)^{-\frac 12}$ to
be the kernel-function for an operator which lies in
$A\times_{\a}G$. In \cite{Rf7} this was achieved simply
by assuming that this kernel-function is in $L^1(G,A)$ for
$a$ and $b$ in a dense subalgebra, much as discussed in
Theorem 5.7 above. But under our present more general hypotheses one
can find examples where $a, b \in {\Cal P}_\a$ but the above 
kernel-function is not in $L^1(G, A)$. This does not mean that
such a kernel-function could not still represent an element of
$A \times_\a G$. But Exel \cite{E2} has shown that in general it does not.
In fact, for the case in which $G$ is Abelian he gives necessary and
sufficient conditions for this to happen. We refer the reader to the
very interesting ``relative continuity'' condition which he shows must
hold, and to his discussion of 
the difficulty of finding a big subspace of mutually relatively
continuous elements.

\section{7.  Square-Integrable Representations}

In this section we study the special case in which the algebra $A$ is the
algebra
$K$ of compact operators.  Then $M(K) = K''$, which very much simplifies
matters.
(For certain considerations a more general setting would involve
$C^*$-algebras $A$
such that $M(A)$ is monotone complete \cite{Pe}.)  The strict topology on
$M(K)$
coincides with the ultra-strong-$*$ operator topology (p.~76 of \cite{La}).
Every
bounded increasing net of self-adjoint elements in $M(K)$ converges in the
strong,
so ultra-strong-$*$ and strict, topologies (lemma $6.1.4$ of \cite{KR1}),
and so for
an action $\alpha$ of a locally compact group $G$, every integrable element is
proper.  That is, ${\Cal P}_{\alpha} = {\Cal M}_{\alpha}$ with the notation
of the
previous sections.

We will show that proper actions are closely related to square-integrable
representations of $G$.  While this is not surprising, it turns out to
provide an
attractive viewpoint on square-integrable representations.

Let $K$ be realized as the algebra, $K(H)$, of compact operators on a
Hilbert space
$H$.  Every automorphism of $K$ is given by conjugation by an element of ${\Cal
U}(H)$, the group of unitary operators on $H$, and this unitary operator is
unique
up to a scalar multiple of modulus $1$.  Thus if $\alpha$ is an action of
$G$ on
$K$, it is given by a projective representation of $G$ on $H$.  For our
purposes
this can be handled most easily 
\cite{Rf3} by passing to the corresponding extension
group.
That is, let
$$
G_{\alpha} = \{(x,u) \in G \times {\Cal U}(H): \alpha_x(a) = uau^{-1}
\text{ for
all } a \in K\}.
$$
Let $T$ denote the group of complex numbers of modulus $1$.  Then we have a
short
exact sequence
$$
0 \to T \to G_{\alpha} \to G \to 0,
$$
where the map {}from $T$ is given by $t \mapsto (e_G,tI_H)$, while the map {}from
$G_{\alpha}$ is given by $(x,u) \mapsto x$.  {}From the topologies on $T$ and
$G$ we
obtain a locally compact topology on $G_{\alpha}$ making the exact sequence of
groups a topological exact sequence.  
The map $(x,u) \mapsto u$ gives an ordinary
unitary
representation of $G_{\alpha}$ on $H$.  The corresponding action on $K$
will be the
pull-back to $G_{\alpha}$ of the action $\alpha$ of $G$.  By passing to
$G_{\alpha}$ we can in this way always assume that $\alpha$ comes {}from an
ordinary
representation.  Because $T$ is compact, it is easily seen that this
passage does
not affect whether the action on $K$ is proper.

Thus {}from now on we always assume that we have an ordinary unitary
representation,
$U$, of $G$, on a Hilbert space $H = H_U$, and that $\alpha$ is the
corresponding
action on $K$.  Our immediate goal is to find necessary and sufficient
conditions
on $U$ such that $\alpha$ is proper.

Suppose now that $a \in {\Cal P}_{\alpha}^+$, $a \ne 0$.  Since $a$ is a
compact
operator and ${\Cal P}_{\alpha}^+$ is a hereditary cone, it follows that
each of
the spectral projections of $a$, and each of their subprojections, is in ${\Cal
P}_{\alpha}^+$.  Thus ${\Cal P}_{\alpha}$ contains enough rank-one
projections to
generate a $C^*$-subalgebra of $K$ containing $a$.  Consequently for many
purposes
we can focus on the rank-one projections in
${\Cal P}_{\alpha}^+$.  Let
$p$ be such, and let $\xi$ be a unit vector spanning the range of $p$.  For any
$\eta,\zeta \in H$ we will write $\<\eta,\zeta\>_K$ for the rank-one operator
determined by $\eta$ and $\zeta$.  For convenience we take the
inner-product on $H$
to be linear in the second variable.  Thus we set
$$
\<\eta,\zeta\>_K\zeta_0 = \eta\<\zeta,\zeta_0\>
$$
for $\zeta_0 \in H$.  Then $p = \<\xi,\xi\>_K$.  Since $p \in {\Cal
P}_{\alpha}$,
there is a constant $k_{\xi}$ such that
$$
\int \lambda(x)\alpha_x(p)dx \le k_{\xi}^2I_H \tag 7.1
$$
for $\lambda \in {\Cal B}$.  Thus for any $\eta \in H$
$$
k_{\xi}^2\|\eta\|^2 \ge \int \lambda(x)\<\alpha_x(p)\eta,\eta\>dx 
= \int \lambda(x) |\<U_x\xi,\eta\>|^2dx.
$$
Because of the definition of ${\Cal B}$ it follows that $x \mapsto
\<U_x\xi,\eta\>$
is in $L^2(G)$.

\demo{\bf 7.2 Notation} For any $\xi,\eta \in H$ we define the
corresponding
coefficient function $c_{\xi\eta}$  by
$$
c_{\xi\eta}(x) = \<U_x\xi,\eta\>.
$$
\enddemo

With this notation, and with $\xi$ now any vector in the range of $p$, we
see {}from
the above that there is a constant $k_{\xi}$ such that
$$
\|c_{\xi\eta}\|_2 \le k_{\xi}\|\eta\| \tag 7.3
$$
for any $\eta \in H$.

Suppose now that $g \in L^1(G) \cap L^2(G)$, and let $U_g$ denote the
integrated
form of $U$.  Then it follows that
$$
\align
|\<U_g\xi,\eta\>| &= |\int {\bar g}(x)c_{\xi\eta}(x)dx| \\
&\le \|g\|_2k_{\xi}\|\eta\|,
\endalign
$$
for every $\eta \in H$.
Consequently,
$$
\|U_g\xi\| \le k_{\xi}\|g\|_2. \tag 7.4
$$

Suppose, conversely, that $\xi \in H$ and that we know that an inequality of
form
$7.4$ holds for all $g \in L^1 \cap L^2(G)$.  Running the above argument
backward,
we obtain $7.3$, and then $7.1$, so that $p \in {\Cal P}_{\alpha}$.  We
have thus
obtained:

\proclaim{7.5 Proposition} Let $U$ be a unitary representation of $G$ on
$H$, with
corresponding action $\alpha$ on $K = K(H)$.  Let $\xi \in H$, and let
$p_{\xi} =
\<\xi,\xi\>_K$.  Then $p_{\xi} \in {\Cal P}_{\alpha}$ if and only if there is a
constant, $k_{\xi}$, such that
$$
\|c_{\xi\eta}\|_2 \le k_{\xi}\|\eta\|
$$
for every $\eta \in H$, or equivalently, such that
$$
\|U_g\xi\| \le k_{\xi}\|g\|_2
$$
for every $g \in L^1 \cap L^2(G)$.
\endproclaim

\definition{7.6 Definition} We will call a vector $\xi$ satisfying these
equivalent
conditions (i\.e\. 7.3 and 7.4) a $U$-{\it bounded} vector.  We will
denote the
set of $U$-bounded vectors by ${\Cal B}_U$.
\enddefinition

This definition is closely related to the definition of bounded elements in Hilbert
algebras \cite{D,Rf1,Pj,Cm}. Compare also with Connes' treatment of
square-integrable representations of foliations beginning on page 573
of \cite{Cn}. (For a recent variation see definition 1.3
of \cite{Bi}.)  It is clear that ${\Cal B}_U$ is a linear
subspace,
possibly not closed, in $H$.  If $\xi \in {\Cal B}_U$ and $x \in G$, then
for $\eta
\in H$
$$
c_{U_x\xi,\eta}(y) = \<U_yU_x\xi,\eta\> = c_{\xi\eta}(yx).
$$
Thus
$$
\|c_{U_x\xi,\eta}\|_2 = \Delta(x)^{-1/2}\|c_{\xi\eta}\|_2. \tag 7.7
$$
Consequently ${\Cal B}_U$ is carried into itself by $U$.

\definition{7.8 Definition} We will say that a unitary representation $U$
of $G$ on
$H$ is {\it square-integrable} if ${\Cal B}_U$ is dense in $H$.
\enddefinition

This is exactly the special case for groups of Combes' definition in $1.7$ of
\cite{Cm} for left Hilbert algebras.

We will see that this definition is equivalent to the more traditional
definitions
in those situations where they have been given.  But conditions $7.3$ and
$7.4$,
which do not seem to have been especially emphasized before, are very
convenient.
In view of our discussion just before $7.1$ of the fact that if $a \in {\Cal
P}_{\alpha}^+$ then all its spectral projections must be in ${\Cal
P}_{\alpha}^+$,
we almost immediately obtain {}from Proposition $7.5$:

\proclaim{7.9 Theorem} Let $U$ be a unitary representation of $G$ on $H$, with
corresponding action $\alpha$ on $K(H)$.  Then $\alpha$ is proper if and
only if  $U$ is square-integrable.
\endproclaim

We now begin to show the relation with the usual definitions of
square-integrable
representations given in the irreducible or cyclic cases
\cite{D, Ro, DM, M, Pj, Ca}.  Let $\xi \in {\Cal
B}_U$.  Define an operator, $T_{\xi}$, {}from $L^1 \cap L^2(G)$ into $H$ by
$$
T_{\xi}(g) = U_g(\xi).
$$
The definition of ${\Cal B}_U$ says that $T_{\xi}$ is bounded for the
$L^2$-norm,
with $\|T_{\xi}\| \le k_{\xi}$.  Furthermore, if we denote by $L$ the left
regular
representation of $G$ on $L^2(G)$, we have
$$
T_{\xi}(L_xg) = U_{L_xg}(\xi) = U_xU_g(\xi) = U_x(T_{\xi}(g)).
$$
Thus $T_{\xi}$ extends to a bounded intertwining operator {}from $L^2(G)$ to $H$,
which we still denote by $T_{\xi}$.  Since $U$ is assumed to be
non-degenerate, it
is clear that the closure of the range of $T_{\xi}$ is exactly the cyclic
subspace
in $H$ generated by $\xi$.  If we now form the polar decomposition of
$T_{\xi}$,
then the partially isometric term will be a unitary intertwining operator
{}from some
closed invariant subspace of $L^2(G)$ onto this cyclic subspace.  (See
VI\.13\.14
of \cite{FD}.)  We thus obtain:

\proclaim{7.10 Proposition} Let $U$ be a unitary representation of $G$, and let
$\xi \in {\Cal B}_U$.  Then the restriction of $U$ to the cyclic subspace
generated by
$\xi$ is unitarily equivalent to a subrepresentation of the left regular
representation of $G$.
\endproclaim

To clarify the situation a bit more, we note the following analogue of
Proposition 5.1, which follows immediately {}from 7.4:

\proclaim{7.11 Proposition} Let $U$ and $V$ be two unitary representations
of $G$,
and let $T$ be a bounded intertwining operator {}from $U$ to $V$.  Then
$$
T({\Cal B}_U) \subseteq {\Cal B}_V.
$$
\endproclaim

\proclaim{7.12 Corollary} Let $U$ be a unitary representation of $G$, and
let $P$
be the projection operator onto a $U$-invariant subspace.  Then $P({\Cal B}_U)
\subseteq {\Cal B}_U$.
\endproclaim

The following proposition is almost immediate {}from Definitions 7.8.

\proclaim{7.13 Proposition} The direct sum of a (possibly infinite) family of
square-integrable unitary representations of $G$ is square-integrable.
\endproclaim

{}From Corollary $7.12$ and the usual process of decomposing a representation
into a
(possibly infinite) direct sum of cyclic representations, we immediately
obtain one
direction of:

\proclaim{7.14 Theorem} The square-integrable representations of $G$ are
exactly
those which are unitarily equivalent to a (possibly infinite) direct sum of
copies
of subrepresentations of the left regular representation of $G$.
\endproclaim

\demo{Proof} We must show the converse.  The crux is to show that the left
regular
representation of $G$ is square-integrable.  Let $\xi \in C_c(G)$.  For any
$g \in
L^1 \cap L^2(G)$ we have
$$
(L_g\xi)(x) = \int g(xy)\xi(y^{-1})dy 
= \int (R_yg)(x)(J{\bar \xi})(y)dy = (R_{J{\bar \xi}}g)(x),
$$
where $(R_yg)(x) = \Delta(y)^{1/2}g(xy)$ so that $R$ is the unitary
right-regular
representation, and $J$ is the Tomita--Takesaki operator \cite{KR} defined by
$(J\eta)(y) = \Delta(y)^{-1/2}\overline{\eta(y^{-1})}$.  Thus
$$
\|L_g\xi\|_2 = \|R_{J{\bar \xi}}g\|_2 \le \|J\xi\|_1\|g\|_2.
$$
(More generally, we see that if $\xi \in L^2(G)$ and $\|J\xi\|_1 < \infty$,
then
$\xi \in {\Cal B}_L$.)  Since $C_c(G)$ is dense in $L^2(G)$, this shows
that $L$ is
square-integrable.  The appearance of $J$ in the above calculation
indicates that
something somewhat interesting is happening.  We will pursue this matter
shortly.

{}From Corollary $7.12$ it follows that every subrepresentation of $L$ is
square-integrable.  The proof is then completed by Proposition $7.13$. \qed
\enddemo

The most common definitions of square-integrable representations just
involves the
condition that $c_{\xi\eta} \in L^2(G)$ for some $\xi,\eta \in H_U$ (and
$c_{\xi\eta} \not\equiv 0$).  Our tiny contribution to this aspect is to
point out
now that by using the basic notion of the graph of an unbounded operator,
we can
avoid explicit use of the theory of unbounded operators and their polar
decomposition when dealing with this condition.  A very similar argument,
involving
an irreducible representation, appears in the appendix of \cite{GMP}.

\proclaim{7.15 Proposition} Let $U$ be a unitary representation of $G$ on
$H$.
Let $\xi,\eta \in H$, and suppose that $c_{\xi\eta} \in L^2(G)$.  Let
$H_{\xi}$
denote the cyclic subspace generated by $\xi$, and replace $\eta$ by its
projection
in this subspace.  Then the restriction of $U$ to the cyclic subspace
$H_{\eta}$
generated by (the new) $\eta$ is unitarily equivalent to a
subrepresentation of $L$.
\endproclaim

\demo{Proof} We can assume that $H = H_{\xi}$.  Note that if $c_{\xi\zeta}
\equiv 0$
for some $\zeta \in H$ then $\zeta = 0$, since $\xi$ is cyclic.  Let
$$
\align
{\Cal D} &= \{\zeta \in H: c_{\xi\zeta} \in L^2(G)\}, \\
\Gamma &= \{(\zeta,c_{\xi\zeta}): \zeta \in {\Cal D}\}.
\endalign
$$
Then $\Gamma$ is a closed subspace of $H \oplus L^2(G)$.  For suppose
$\{\zeta_n\}$ is a sequence in ${\Cal D}$ such that $\zeta_n$ converges to
$\zeta
\in H$ and $c_{\xi\zeta_n}$ converges to $h \in L^2(G)$.  Then {}from the
definition of $c_{\xi\zeta_n}$ we see that $\{c_{\xi\zeta_n}\}$ converges
uniformly
pointwise to $c_{\xi\zeta}$.  Thus $c_{\xi\zeta} = h$, so $c_{\xi\zeta} \in
L^2(G)$.

It is easily checked that $\Gamma$ is a $U \oplus L$ invariant subspace of $H
\oplus L^2(G)$.  Consider the operator $Q$ with domain $H$ defined by
$$
\zeta \mapsto (\zeta,0) \mapsto (\zeta',c_{\xi\zeta'}) \mapsto c_{\xi\zeta'}
$$
where the second arrow is the orthogonal projection {}from $H \oplus
L^2(G)$ onto
$\Gamma$.  Then $Q$ is clearly a norm-decreasing operator which intertwines
$U$ and
$L$.  If $Q\zeta = 0$ then $\zeta$ is clearly orthogonal to ${\Cal D}$.
Thus $Q$
is injective on the closure, ${\bar {\Cal D}}$, of ${\Cal D}$, which is a
$U$-invariant subspace.  {}From the polar decomposition of the bounded
operator $Q$
we obtain an isometric intertwining operator {}from ${\bar {\Cal D}}$ to
$L^2(G)$.
Clearly $H_{\eta} \subseteq {\bar {\Cal D}}$. \qed
\enddemo

We remark that if $U$ is irreducible, then $Q$ is already a multiple of an
isometry
{}from ${\bar {\Cal D}}$.  If $G$ is unimodular, then $\|c_{\xi\eta}\|_2 =
\|c_{\eta\xi}\|_2$ in all $\xi,\eta \in H$.  {}From this it follows easily that
${\Cal B}_U = H$ in this case (for irreducible $U$).

We also remark that if $N$ denotes the operator $\zeta \mapsto \zeta'$ for
$\zeta'$
as above, then $N$ is a norm-decreasing intertwining operator on $H$, and
that
$Q\zeta = c_{\xi N\zeta}$.  But
$$
c_{\xi N\zeta}(x) = \<U_x\xi,N\zeta\> = \<U_xN\xi,\zeta\> = c_{N\xi,\zeta},
$$
so $N\xi \in {\Cal B}_U$ in view of the properties of $Q$.

If $U$ is irreducible, then $N$ must be a multiple of an isometry {}from
$H$.  In
particular, every element of $H$ would be in the range of $N$, so that
$c_{\xi\zeta} \in L^2(G)$ for all $\zeta \in H$. When this is combined
with part {\it iii} of the restatement in \cite{BT} of theorem 3
of \cite{DM}, this says that if $\xi$ is ``admissible'' \cite{BT},
then $\xi \in {\Cal B}_U$.

The next proposition ties the situation a bit more closely to the
discussion of the
previous sections.  Its proof is an immediate application of the definitions.

\proclaim{7.16 Proposition} Let $U$ be a unitary representation of $G$ on
$H$, with
corresponding action $\alpha$ on $K$, and let $\xi,\eta \in H$.  Then
$\<\eta,\xi\>_K \in {\Cal Q}_{\alpha}(= {\Cal N}_{\alpha})$ 
iff $\xi \in {\Cal B}_U$.
\endproclaim

The left regular representation of $G$ comes {}from the action of $G$ on
$C_{\infty}(G)$ by left translation, which is proper, together with the
invariant
unbounded (for $G$ not compact) 
trace on $C_{\infty}(G)$ consisting of Haar measure.  This
suggests that
perhaps we obtain square-integrable representations {}from other proper
actions and
invariant traces.  But the occurrence of the operator $J$ in the proof of
Theorem
$3.14$ should warn us of possible difficulties.  On the other hand, because
traces
are ``measure-theoretic'' we will see that we can deal with integrable
actions --
the full force of being proper is not important here.

Let $\alpha$ be an action of $G$ on a $C^*$-algebra $A$.  We recall
\cite{KR2,Pe2}
that a trace on $A$, possibly unbounded, is a function $\tau$ {}from $A^+$ to
$[0,\infty]$ with the expected properties.  The correct set-up for us here
appears to
involve
the following definition.  (See \cite{DM}.)

\definition{7.17 Definition} Let $\tau$ be a trace on $A$.  Then $\tau$ is
said to
be $\Delta$-semi-invariant for the action $\alpha$ if
$$
\tau(\alpha_x(a)) = \Delta(x)^{-1}\tau(a)
$$
for all $a \in A^+$ and $x\in G$.
\enddefinition

Much as earlier, we set ${\Cal M}_{\tau}^+ = \{a \in A^+: \tau(a) < \infty\}$,
${\Cal N}_{\tau} = \{a \in A: a^*a \in {\Cal M}_{\tau}^+\}$, and 
${\Cal M}_{\tau} = \text{span} {\Cal M}_{\tau}^+$. Then
${\Cal M}_{\tau} = {\Cal N}_{\tau}^2$, and ${\Cal M}_{\tau} \cap A^+ = {\Cal
M}_{\tau}^+$, and $\tau$ extends to ${\Cal M}_{\tau}$. Notice that if
$\tau$ is $\Delta$-semi-invariant for $\a$, then ${\Cal M}_\tau^+$
is carried into itself by $\a$, and similarly for ${\Cal N}_\tau$
and ${\Cal M}_\tau$.

We recall (proposition $6.1.3$ of \cite{Pe2}) that if $\tau$ is a
lower semi-continuous trace (or weight) on $A$, then the GNS construction
works to
produce a non-degenerate $*$-representation of $A$.  We denote its Hilbert
space by
$H_{\tau}$, but do not use specific notation for the representation (i\.e\.
we use
module notation).  Each element $a \in {\Cal N}_{\tau}$ determines an
element of
$H_{\tau}$, but our notation will not distinguish between elements of ${\Cal
N}_{\tau}$ and their images in $H_{\tau}$. The first parts of the
following theorem are basically well-known.

\proclaim{7.18 Theorem} Let $\alpha$ be an action of $G$ on a $C^*$-algebra
$A$,
and let $\tau$ be a $\Delta$-semi-invariant lower semi-continuous trace on $A$,
with {\rm GNS} representation on $H_{\tau}$.  Define a unitary
representation $U$
of $G$ on $H_{\tau}$ by
$$
U_x(a) = \Delta(x)^{1/2}\alpha_x(a)
$$
for $a \in {\Cal N}_{\tau}$. Then $U$ is strongly continuous.
Furthermore, 
every element of ${\Cal M}_{\alpha} \cap {\Cal M}_{\tau}$ is $U$-bounded.  If
$\alpha$ is integrable, then $U$ is square-integrable.
\endproclaim

\demo{Proof} We remark that ${\Cal N}_{\tau}$ may be very small, even just
$\{0\}$.  But because $\tau$ is a trace, ${\Cal M}_{\tau}$ is a two-sided
ideal in
$A$, as is then ${\Cal N}_{\tau}$.  (See $6.2.1$ of \cite{Pe2}.)

By the semi-continuity of $\tau$ the image of ${\Cal M}_{\tau}$ in ${\Cal
H}_{\tau}$ is dense (see $7.4.1$ of \cite{D}), and the representation of
$A$ on
${\Cal H}_{\tau}$ is non-degenerate.  For $a \in {\Cal M}_{\tau}$, $b \in {\Cal
M}_{\tau}^+$, and $x \in G$ we have
$$
\<U_x(a),b\>_{\tau} = \tau(\alpha_x(a^*)b) = \tau(b^{1/2}\alpha_x(a^*)b^{1/2})
$$
by proposition 5.2.2 of \cite{Pe2}.  
But $c \mapsto \tau(b^{1/2}cb^{1/2})$
is a
positive linear functional defined on all of $A$, so continuous.  {}From this
and the
fact that ${\Cal M}_{\tau}$ is the span of ${\Cal M}_{\tau}^+$ we see that
$U$ is
weakly continuous.  Thus $U$ is strongly continuous since it is unitary.

Suppose now that $a \in {\Cal M}_{\alpha}^+ \cap {\Cal M}_{\tau}$ and that
$g \in
L^1 \cap L^2(G)$.  Assume further that $g$ has compact support.  Since $U$ is
strongly continuous, the integrated form, $U_g$, is defined.  Then
$$
\align
\|U_ga\|^2 &= \<a,U_{g^**g}a\>_{\tau} 
= \tau(a^{1/2}\alpha_{\Delta^{1/2}(g^**g)}(a)a^{1/2}) \\
&\le \tau(a)\|\alpha_{\Delta^{1/2}(g^**g)}(a)\| 
\le \tau(a)k_a\|\Delta^{1/2}(g^**g)\|_{\infty},
\endalign
$$
since $a \in {\Cal M}_{\alpha}^+$.  But for each $x \in G$ we have
$$
\Delta^{1/2}(x)|(g^**g)(x)| \le \Delta^{1/2}(x) \int |{\bar g}(y)g(yx)|dy 
\le \Delta^{1/2}(x)\|g\|_2\|g(\cdot x)\|_2 = \|g\|_2^2.
$$
Thus $\|\Delta^{1/2}(g^**g)\|_{\infty} \le \|g\|_2^2$.  Putting this
together, we
obtain
$$
\|U_ga\|^2 \le \tau(a)k_a\|g\|_2^2.
$$
Since any $g \in L^1 \cap L^2(G)$ can be approximated by ones of compact
support,
simultaneously in $L^1$ and $L^2$ norm, this inequality holds for all $g
\in L^1
\cap L^2(G)$.  This says that $a$ as vector in ${\Cal H}_{\tau}$ is
$U$-bounded.
Thus we see that every element of ${\Cal M}_{\alpha}^+ \cap {\Cal M}_{\tau}$ is
$U$-bounded.  Since ${\Cal M}_{\tau}$ is the span of ${\Cal M}_{\tau}^+$, and
similarly for ${\Cal M}_{\alpha}$, it follows that every element of ${\Cal
M}_{\alpha} \cap {\Cal M}_{\tau}$ is $U$-bounded.  

Suppose now that $\alpha$ is
integrable.  To show that $U$ is square-integrable it suffices to show that
${\Cal
M}_{\alpha} \cap {\Cal M}_{\tau}$ is dense in ${\Cal H}_{\tau}$.  Since ${\Cal
M}_{\tau}$ is an ideal and ${\Cal M}_{\alpha}$ is hereditary, we see that
${\Cal
M}_{\alpha}{\Cal M}_{\tau}{\Cal M}_{\alpha} \subset {\Cal M}_{\alpha} \cap
{\Cal
M}_{\tau}$.  We show that in ${\Cal H}_{\tau}$ every element of ${\Cal
M}_{\tau}$
can be approximated by elements of ${\Cal M}_{\alpha}{\Cal M}_{\tau}{\Cal
M}_{\alpha}$.  Since, as noted above, ${\Cal M}_{\tau}$ is dense in ${\Cal
H}_{\tau}$, this will conclude the proof.
\enddemo

Let $I$ denote the norm-closure of ${\Cal M}_{\tau}$ in $A$, so that $I$ is an
$\alpha$-invariant ideal in $A$.  Since ${\Cal M}_{\alpha}$ is assumed dense in
$A$, it follows that ${\Cal M}_{\alpha}I{\Cal M}_{\alpha}$ is dense in the
$C^*$-algebra $I$.  But ${\Cal M}_{\alpha}I{\Cal M}_{\alpha} \subset {\Cal
M}_{\alpha}$ since ${\Cal M}_{\alpha}$ is hereditary.  Thus ${\Cal
M}_{\alpha} \cap
I$ is dense in $I$.  Although we don't need it here, we note that we have
essentially proven the following analogue of Proposition 5.5:

\proclaim{7.19 Proposition} Let $\alpha$ be an integrable action of $G$ on
$A$, and
let $I$ be an $\alpha$-invariant ideal in $A$.  Then the action $\alpha$ of
$G$ on
$I$ is integrable.
\endproclaim

We continue with the proof of Theorem $7.18$.  Pick a positive approximate
identity
of norm $1$ for $I$.  Since ${\Cal M}_{\alpha} \cap I$ is a dense
$*$-subalgebra of
$I$, we can approximate the approximate identity by elements of ${\Cal
M}_{\alpha}
\cap I$ to obtain a self-adjoint approximate identity of norm $1$ consisting of
elements of ${\Cal M}_{\alpha} \cap I$.  We can then square this approximate
identity to obtain one which is positive. We denote the resulting
approximate identity in ${\Cal M}_{\alpha} \cap I$ by $\{e_{\lambda}\}$.

Let $b \in {\Cal M}_{\tau}^+$.  Then $e_{\lambda}be_{\lambda} \in {\Cal
M}_{\alpha}
\cap {\Cal M}_{\tau}$.  We show that $\{e_{\lambda}be_{\lambda}\}$
converges to $b$
in ${\Cal H}_{\tau}$.  Now, using heavily that $\tau$ is tracial, we have
$$
\|b-e_{\lambda}be_{\lambda}\|_{\tau}^2 = \tau(b^2-2be_{\lambda}be_{\lambda} +
be_{\lambda}^2be_{\lambda}^2) 
= \tau(b^2) - \tau(b^{1/2}(2e_{\lambda}be_{\lambda} -
e_{\lambda}^2be_{\lambda}^2)b^{1/2}).
$$
But $b^{1/2} \in {\Cal N}_{\tau}$ and so $a \mapsto \tau(b^{1/2}ab^{1/2})$ is a
positive linear functional defined everywhere on $A$, and so continuous.  Since
$2e_{\lambda}be_{\lambda} - e_{\lambda}^2be_{\lambda}^2$ converges to $b$
in norm,
we see that $e_{\lambda}be_{\lambda}$ does indeed converge to $b$ in ${\Cal
H}_{\tau}$. \qed

It is not at all clear to me how much of the above can be done if $\tau$ is
only a
weight instead of a trace.  As mentioned before Question 6.9, 
even the question as to whether the unitary
representation $U$ is strongly continuous seems quite delicate. 

Let $\alpha$ now be the proper action of $G$ on $A = C_{\infty}(G)$ by {\it
right}
translation, so $(\alpha_x(f))(y) = f(yx)$.  Let $\tau$ be the trace on $A$
defined
by {\it left} Haar measure.  Then $\pi$ is $\Delta$-semi-invariant for
$\alpha$.
Thus we obtain what we already know:

\proclaim{7.20 Corollary} The right regular representation of $G$ on
$L^2(G)$ is
square-integrable.
\endproclaim

But we know that the left-regular representation too is square-integrable.
(It is
equivalent to the right regular representation.)  Even more, the left action
of a subgroup of $G$ on $L^2(G)$ should be square-integrable. 
We can relate this to Theorem
$7.18$ as follows.  Let $\alpha$ be an action of $G$ on a $C^*$-algebra
$A$, and
let $\tau$ be a trace on $A$ which we now suppose to be actually
$\alpha$-invariant.  Suppose that $d$ is an unbounded positive invertible
operator
affiliated with $A$ in the sense of Woronowicz \cite{Wo}.  (See also Baaj
\cite{Ba}.)  For our purposes this means that we have a morphism, say $\theta$,
{}from $D = C_{\infty}({\Bbb R})$ to $A$ (that is, a $*$-homomorphism {}from
$D$ into
$M(A)$ such that $\theta(D)A$ is dense in $A$) together with a strictly
positive
$\delta \in C({\Bbb R})$ (where $C({\Bbb R})$ denotes the algebra of possibly
unbounded continuous function on ${\Bbb R}$) acting by pointwise
multiplication as
an unbounded operator on $C_{\infty}({\Bbb R})$ with domain $D_0 =
C_c({\Bbb R})$.
Then $d$ is defined to be the closure of the operator on $A$ with domain
$\theta(D_0)A$ defined by
$$
d(\theta(\varphi)a) = (\theta(\delta\varphi))a
$$
for $\varphi \in D_0$, $a \in A$.  Let $C$ denote the range of $\theta$,
and set
$C_0 = \theta(D_0)$.  For our present purposes we require that $d$ be
central, that
is, that $C \subseteq ZM(A)$, where $ZM(A)$ denotes the center of $M(A)$.  
Now $\alpha$ carries $ZM(A)$ into itself,
and $d$
can be represented by an unbounded continuous function on the maximal ideal
space
of the center.  {}From this point of view it is clear what we mean by
$\alpha_x(d)$
for $x \in G$.

\definition{7.21 Definition} We say that a central positive operator $d$
affiliated with
$A$, via the morphism $\theta$ {}from $C_{\infty}({\Bbb R})$ to $A$, is
$\Delta$-{\it
semi-invariant} for $\alpha$ if the range of $\theta$ is carried into itself by
$\alpha_x$ and
$$
\alpha_x(d) = \Delta(x)d
$$
for all $x \in G$.
\enddefinition

As one example we have:

\proclaim{7.22 Proposition} Let $\a$ be a proper action of $G$ on a
locally compact space $M$, and so on $A = C_\infty(M)$. Assume that
$M/\a$ is paracompact. Then there exists a central positive
invertible operator affiliated with $A$ which is $\Delta$-semi-invariant
for $\a$.
\endproclaim

\demo{Proof} We imitate the construction of ``Bruhat approximate
cross-sections''. For each $\a$-orbit choose an element of
$C_c(M)^+$ which is not everywhere 0 on that orbit. The images in
$M/\a$ of the sets where these functions are non-zero form an open
cover of $M/\a$. By paracompactness there is a locally finite subcover.
Let $b$ denote the sum of the functions for this subcover. So $b$ is
a continuous positive function, with the property that its support
meets the preimage in $M$ of any compact subset of $M/\a$ in a 
compact set, and it is not everywhere 0 on any orbit. Define a
function $h$ on $M$ by
$$
  h(m) = \int_G b(\a_y^{-1}(m))\Delta^{-1}(y)\ dy   .
$$
The integrand has compact support for each $m$, so $h$ is
well-defined. {}From the properties of $b$ it is clear that $h$ is 
positive, continuous, and nowhere 0, so invertible. Furthermore,
for $x \in G$ we have
$$
(\a_x(h))(m) = \int b(\a_y^{-1}(\a_x^{-1}(m)))\Delta^{-1}(y)\ dy 
= \Delta(x)h(m)  .
$$
Thus $h$ is $\Delta$-semi-invariant. (The requirement on domains
is easily checked.)    \qed
\enddemo

It is natural to wonder whether there are interesting extensions of this
construction for non-commutative $A$'s. 

We now continue the general development.

For $c \in ZM(A)^+$ define $\tau_c$ on $A^+$ by
$$
\tau_c(a) = \tau(ca) = \tau(a^{1/2}ca^{1/2}).
$$
It is clear that $\tau_c$ is a trace on $A$.  Note that if $a \in {\Cal
M}_{\tau}$
then $ca \in {\Cal M}_{\tau}$ since $ca = a^{1/2}ca^{1/2} \le \|c\|a$ and
${\Cal
M}_{\tau}$ is hereditary.  Thus ${\Cal M}_{\tau_c} \supseteq {\Cal
M}_{\tau}$.  It
is easily seen that since $\tau$ is lower semi-continuous, so is $\tau_c$.

Now let $d$ be a central 
positive operator affiliated with $A$, by means of the morphism
$\theta$ {}from $D$ to $A$ and $\delta \in C({\Bbb R})$.  Analogously to what
we did
in the previous sections let ${\Cal B}(C_0) = \{c \in C_0: 0 \le c \le
1\}$, with
its usual upward directed order.  Then $\{\tau_{dc}: c \in {\Cal B}(C_0)\}$
is an
increasing net of lower semi-continuous traces on $A$, and so we can define
$\tau_d$ to be their upper bound.  Thus $\tau_d$ is a lower semi-continuous
trace
on $A$.  For $a \in A^+$ and $c_0 \in C_0^+$ we have
$$
\tau_d(c_0a) = \lim_c \tau(dcc_0a) = \lim_c \tau(a^{1/2}dcc_0a^{1/2}) 
= \tau((dc_0)a)
$$
where $c$ ranges over ${\Cal B}(C_0)$.  In particular, $C_0{\Cal M}_{\tau}
\subseteq {\Cal M}_{\tau_d}$.

We are assuming that $d$ is invertible.  It is then easy to see that $\tau$
comes
{}from $\tau_d$ by the above procedure using $d^{-1}$.  That is, $\tau =
(\tau_d)_{d^{-1}}$.  In particular, $C_0{\Cal M}_{\tau_d} \subseteq {\Cal
M}_{\tau}$.  Now $C_0 = C_0^3$, so
$$
C_0{\Cal M}_{\tau} = C_0^3{\Cal M}_{\tau} \subseteq C_0^2{\Cal M}_{\tau_d}
\subseteq C_0{\Cal M}_{\tau}.
$$
Thus $C_0{\Cal M}_{\tau} = C_0{\Cal M}_{\tau_d}$.  By a calculation which is
computationally simpler than that near the end of the proof of Theorem
$7.18$ one
sees that $C_0{\Cal M}_{\tau}$ is dense in ${\Cal H}_{\tau}$.  One must
just notice
that the fact that $C_0$ is in $M(A)$ rather than $A$ causes no
difficulties.  {}From
what we have seen, $C_0{\Cal M}_{\tau}$ is also dense in ${\Cal
H}_{\tau_d}$.  Let
us define an operator, $T$, {}from ${\Cal H}_{\tau_d}$ to ${\Cal H}_{\tau}$
by first
defining it on $C_0{\Cal M}_{\tau}$ by
$$
T(ca) = (d^{1/2}c)a.
$$
One checks that this is well-defined as follows.  Given $\sum c_ia_i = \sum
c'_ja'_j$, choose $c \in C_0$ such that $cc_i = c_i$ and $cc'_j = c'_j$ for all
$i,j$.  Then
$$
T(\sum c_ia_i) = (d^{1/2}c)\sum c_ia_i = T(\sum c'_ja'_j).
$$
Then
$$
\<Tc_1a_1, Tc_2a_2\>_{\tau} = \tau((dc_2^*c_1)a_2^*a_1) 
= \<c_1a_1,c_2a_2\>_{\tau_d}.
$$
Thus on its domain $T$ is isometric.  But its domain and range are dense in
${\Cal
H}_{\tau_d}$ and ${\Cal H}_{\tau}$ respectively.  So $T$ extends to a unitary
operator between them.

Now suppose that $\alpha$ is an action of $G$ on $A$, that $\tau$ is
$\alpha$-invariant, and that $d$ is $\Delta$-semi-invariant for $\alpha$.  In
particular, $\alpha$ carries $C_0$ into itself and is an automorphism of the
directed set ${\Cal B}(C_0)$.  Then for $c \in {\Cal B}(C_0)$ and $a \in
A^+$ we
have
$$
\align
\tau_{dc}(\alpha_x(a)) &= \tau(dc\alpha_x(a)) 
= \tau(\alpha_x((\alpha_{x^{-1}}(dc))a)) \\
&= \tau(\alpha_{x^{-1}}(dc)a) 
= \Delta(x)^{-1}\tau(d\alpha_{x^{-1}}(c)a) 
= \Delta(x)^{-1}\tau_{d\alpha_{x^{-1}}(c)}(a).
\endalign
$$
On taking the limit over ${\Cal B}(C_0)$ we obtain
$$
\tau_d(\alpha_x(a)) = \Delta(x^{-1})\tau_d(a),
$$
that is, $\tau_d$ is $\Delta$-semi-invariant for $\alpha$.

We can now apply Theorem $7.18$ to conclude that the unitary representation
$V$ on
${\Cal H}_{\tau_d}$ coming {}from $\alpha$ is square-integrable if $\alpha$ is
integrable.  Of course $V$ is defined by
$$
V_x(a) = \Delta(x)^{1/2}\alpha_x(a).
$$
At the same time we have the unitary representation $U$ on ${\Cal H}_{\tau}$
defined by
$$
U_x(a) = \alpha_x(a).
$$
But consider the unitary operator $T$ defined several paragraphs ago.  For
$c \in
C_0$ and $a \in {\Cal M}_{\tau}$ we have
$$
\align
U_x(T(ca)) &= U_x((d^{1/2}c)a) = \alpha_x((d^{1/2}c)a) 
= \alpha_x(d^{1/2}c)\alpha_x(a) \\
&= \Delta(x)^{1/2}d^{1/2}\alpha_x(c)\alpha_x(a) 
= T(\Delta(x)^{1/2}\alpha_x(ca)) = T(V_x(ca)).
\endalign
$$
Thus $T$ is an intertwining operator, and the two representations are
equivalent.
We have thus demonstrated:

\proclaim{7.23 Theorem} Let $\alpha$ be an action of $G$ on a $C^*$-algebra
$A$,
and let $\tau$ be an $\alpha$-invariant lower semi-continuous trace on $A$.
Let
$U$ be the corresponding unitary representation of $G$ on ${\Cal
H}_{\tau}$.  If
$\alpha$ is integrable, and if there is a central positive invertible operator
affiliated to
$A$ which is $\Delta$-semi-invariant, then $U$ is square-integrable.
\endproclaim

We will see a reflection of this theorem in the next section. Upon
applying Proposition 7.22 we obtain:

\proclaim{7.24 Corollary} Let $\a$ be a proper action of $G$ on a 
locally compact space $M$ such that $M/\a$ is paracompact. For
every positive $\a$-invariant Radon measure $\mu$ on $M$ the
corresponding unitary representation of $G$ on $L^2(M, \mu)$ is
square-integrable.
\endproclaim

As one (unsurprising) 
application of the earlier Theorem $7.18$ we can consider the canonical
trace on the algebra of compact operators, whose GNS Hilbert space is the
space of
Hilbert--Schmidt operators.

\proclaim{7.25 Corollary} Let $G$ be a unimodular group and let $U$ be a
square-integrable representation of $G$ on $H$.  Let $\alpha$ be the
corresponding action on $K(H)$, which is integrable.  Then the corresponding
unitary representation on the space of Hilbert--Schmidt operators is
square-integrable.
\endproclaim

\section{8.  The Orthogonality Relations}

In this section we examine what the orthogonality relations for
square-integrable
representations look like {}from the vantage point of the previous section.

Let $U$ be a representation of $G$ on $H$, and let $\alpha$ be the
corresponding
action on $K = K(H)$.  Then $M(K)$ consists of the bounded operators
on $H$, and
$M(K)^{\alpha}$
is exactly the algebra of intertwining operators for $U$.  Let $\xi$ and
$\eta \in
{\Cal B}_U$.  {}From Proposition $7.5$ and polarization it follows rapidly that
$\<\xi,\eta\>_K \in {\Cal M}_{\alpha}$, and so the integral
$$
\int \alpha_x(\<\xi,\eta\>_K)dx \tag 8.1
$$
converges in the strong operator topology to an operator in
$M(K)^{\alpha}$.  {}From
the vantage point of this paper, the orthogonality relations are concerned with
identifying to some extent this intertwining operator.  The reason is that
for any
$\zeta,\omega \in H$ we have
$$
\aligned
\<\int \alpha_x(\<\xi,\eta\>_K)dx\ \zeta,\omega\> &= \int
\overline{\<U_x\eta,\zeta\>} \<U_x\xi,\omega\>dx \\
&= \<c_{\eta\zeta},c_{\xi\omega}\>_{L^2(G)},
\endaligned \tag 8.2
$$
so that any answer will say something about the inner-product of the
coefficient
functions.  It is quite clear that if $\xi$ and $\eta$ come {}from two
subrepresentations which are disjoint (have no non-zero intertwining operators)
then $8.1$ must be $0$ and so $c_{\xi\omega}$ and $c_{\eta\zeta}$ must be
orthogonal.  Since any two representations can be viewed as
subrepresentations of
their direct sum, we obtain:

\proclaim{8.3 Proposition} (The ``first orthogonality relation''.)  Let $U$
and $V$
be representations of $G$ on $H_U$ and $H_V$.  Suppose that $U$ and $V$ are
disjoint.  Then for any $\xi,\omega \in H_U$ with $\xi$ $U$-bounded, 
and for any
$\eta,\zeta \in H_V$ with $\eta$ $V$-bounded, the coefficient functions
$c_{\xi\omega}$ and $c_{\eta\zeta}$ are orthogonal in $L^2(G)$.
\endproclaim

So we concentrate on the ``second orthogonality relation''.  We introduced
after
Theorem $7.9$ the bounded operators $T_{\xi}$ for $\xi \in {\Cal B}_U$.  Let us
calculate $T_{\xi}^*$.  For $\eta \in H$ and $g \in L^1 \cap L^2(G)$ we have
$$
\<g,T_{\xi}^*\eta\> = \<L_g\xi,\eta\> = \<g,c_{\xi\eta}\>.
$$
Thus
$$
T_{\xi}^*\eta = c_{\xi\eta}.
$$
For $\xi,\eta \in {\Cal B}_U$ and $\zeta,\omega \in H$ it follows that
$$
\<T_{\xi}T_{\eta}^*\zeta,\omega\> = \<c_{\eta\zeta},c_{\xi\omega}\>.
$$
{}From $8.2$ we then obtain:

\proclaim{8.4 Proposition} Let $U$ be a unitary representation of $G$.  For
$\xi,\eta \in {\Cal B}_U$ we have
$$
\int \alpha_x(\<\xi,\eta\>_K)dx = T_{\xi}T_{\eta}^*.
$$
\endproclaim

We now examine what this says for the left regular representation, $L$, of
$G$.
Let $h \in {\Cal B}_L$.  {}From the calculation in the proof of Theorem
$7.14$, but
now with $g \in C_c(G)$, we see that, at the level of functions,
$$
T_hg = R_{J{\bar h}}g.
$$
In interpreting this, note that $J$ is an isometry, so that $J{\bar h} \in
L^2(G)$.  The operator $R_{J{\bar h}}$, which is initially defined only on,
say,
$C_c(G)$, is bounded, by our assumption on $h$, and so extends to a bounded
operator on all of $L^2(G)$.  With this understanding, we have
$$
T_h = R_{J{\bar h}}.
$$
The second orthogonality relation for the left regular representation can
then be
considered to be the following statement:

\proclaim{8.5 Theorem} For $f,g \in {\Cal B}_L$ we have
$$
\int \alpha_x(\<f,g\>_K)dx = R_{J{\bar f}}R_{J{\bar g}}^*.
$$
\endproclaim

We remark that this is closely related to the result in example $2.1$ of
\cite{Rf7} (where the $\rho$ there is not unitary).

We now relate this to Exel's example of section 13 of \cite{E2}. Suppose
that $G$ is Abelian. We conjugate $L$ by the Fourier transform so that
it acts on $L^2(\hat G)$, by pointwise multiplication of characters. Upon
applying the Plancherel theorem, we see that equation 7.4 becomes

$$
         \|{\hat g}\xi\|_2 \leq k_\xi\|\hat g\|_2   
$$
for all $g \in L^1\cap L^2(G)$ and $\xi \in L^2(\hat G)$. This will hold
exactly if  $\xi \in L^\infty(\hat G)$. Thus under this picture
${\Cal B}_L = L^\infty\cap L^2(\hat G)$. {}From the fact that $G$ is Abelian
it is easily seen that  $R_{J\bar f} = L_f$, so that the relation in Theorem
8.5 reads
$$
      \int \alpha_x(\<f,g\>_K)dx = L_{f\star g^*}   .
$$
But on $L^2(\hat G)$ the operator $L_{f\star g^*}$ is just pointwise
multiplication by the Fourier transform of $f\star g^*$. If we change 
notation so that now $f, g \in L^\infty\cap L^2(\hat G) = {\Cal B}_L$, and
if we let $M_{f\bar g}$ denote the operator of pointwise multiplication by
$f\bar g$, we obtain
$$
          \int \alpha_x(\<f,g\>_K)dx = M_{f\bar g}   .
$$
Let $f, g \in L^\infty\cap L^2(\hat G)$, and choose any $h  \in L^\infty\cap L^2(\hat G)$
such that $\|h\|_2 = 1$. Then
$$
     \<f, g\>_K = \<f, h\>_K\<h, g\>_K   ,
$$
and so {}from Proposition 7.5 and the fact that ${\Cal P}_\alpha$  is an algebra
we see that the ``big fixed-point-algebra'' will contain $M_{f\bar g}$. In particular,
if $\hat G$ is compact, the big fixed-point algebra is exactly $L^\infty(\hat G)$. This
is Exel's example. Certainly $L^\infty(\hat G)$ is too big. As Exel indicates,
it is well-known that $K(L^2(G))\times_\alpha G$ is isomorphic to
$K(L^2(G))\otimes C_\infty(\hat G)$, which has $\hat G$ as primitive ideal
space. Since the primitive ideal spaces of strongly Morita equivalent
algebras are homeomorphic \cite{Rf3}, it is impossible for $L^\infty(\hat G)$ to
be strongly Morita equivalent to an ideal in $K\times_\alpha G$ (unless $G$ is
finite). We refer the reader to \cite{E2} for substantial further exploration
of this situation. 

We wish to consider next the case in which $U$ is irreducible.  But we
first make
some observations about the general case which will be useful for that
purpose.
For any $\xi,\eta \in {\Cal B}_U$, the operator $T_{\xi}^*T_{\eta}$ on $L^2(G)$
intertwines $L$.  For any $g \in L^1 \cap L^2(G)$ we have
$$
T_{\xi}^*T_{\eta}g = T_{\xi}^*U_g\eta = L_gT_{\xi}^*\eta = g * c_{\xi\eta}.
$$
Note that this implies that $c_{\xi\eta}$ is in the closure of the range of
$T_{\xi}^*T_{\eta}$.  Note next that $c_{\xi\eta}(x^{-1}) = {\bar
c}_{\eta\xi}(x)$,
so that $c_{\xi\eta}$ is in $L^2(G,\Delta^{-1}dx)$ as well as in $L^2(G)$,
where by
$\Delta^{-1}dx$ we denote right Haar measure.  Now for any $\varphi \in
L^2(G)$ and
$\psi \in L^2(G,\Delta^{-1}dx)$ we have
$$
|\varphi*\psi(x)| = |\int \varphi(y)\psi(y^{-1}x)dy| \le
\|\varphi\|_2\|\psi(\cdot^{-1})\|_2.
$$
This says that convolution is well defined and jointly continuous {}from $L^2(G)
\times L^2(G,\Delta^{-1}dx)$ to $L^{\infty}(G)$.  But if $\varphi,\psi \in
C_c(G)$,
then $\varphi*\psi \in C_c(G) \subseteq C_{\infty}(G)$.  Since $C_c(G)$ is
dense in
$L^2(G)$, it follows that $\varphi *\psi \in C_{\infty}(G)$ for all
$\varphi \in
L^2(G)$ and $\psi \in L^2(G,\Delta^{-1}dx)$.

In particular, for $\xi,\eta \in {\Cal B}_U$ and for any $\varphi \in
L^2(G)$, the
function $\varphi * c_{\xi\eta}$ is continuous.  But we saw above that
$g*c_{\xi\eta} = T_{\xi}^*T_{\eta}g$ for $g \in L^1 \cap L^2(G)$.  Let
$\{g_n\}$ be
a sequence in $L^1 \cap L^2(G)$ which converges to $\varphi$.  As seen
above, $g_n
* c_{\xi\eta}$ converges uniformly to $\varphi * c_{\xi\eta}$, but it also
converges in $L^2$-norm to $T_{\xi}^*T_{\eta}\varphi$.  (The main concern
here is
the fact that $c_{\xi\xi}$ need not be in $L^1(G)$.)  We thus obtain:

\proclaim{8.6 Proposition} Let $U$ be a unitary representation of $G$.  For any
$\xi,\eta \in {\Cal B}_U$ and any $\varphi \in L^2(G)$,
$$
T_{\xi}^*T_{\eta}\varphi = \varphi * c_{\xi\eta},
$$
and $\varphi * c_{\xi\eta} \in C_{\infty}(G) \cap L^2(G)$.  Thus the range of
$T_{\xi}^*T_{\eta}$ consists entirely of functions in $C_{\infty}(G)$.
\endproclaim

We now consider the case in which $U$ is an irreducible representation of
$G$.  In
this case, since ${\Cal B}_U$ is an invariant subspace, it must be dense in
$H$
as soon as it contains one non-zero vector, which we will assume.  Let
$\xi,\eta
\in {\Cal B}_U$.  Since $T_{\xi}T_{\eta}^*$ is an intertwining operator, it
must be
a scalar multiple of the identity operator.  We denote the scalar by
$\gamma(\xi,\eta)$, so that
$$
    T_\xi T_\eta ^* = \gamma(\xi, \eta)I_H    .
$$
We wish to obtain a more revealing expression for $\gamma$.  We
follow
the general outline of the treatment given in \cite{Ca}, but our details
are more
elementary because of our use of ${\Cal B}_U$.

By suitably normalizing $\xi$, we can arrange that $T_{\xi}^*$ is an isometry.
Then $T_{\xi}^*T_{\xi}$ is a projection operator on $L^2(G)$, intertwining
$L$.
The restriction of $L$ to the range of $T_{\xi}^*T_{\xi}$ is a
subrepresentation of
$L$ which is unitarily equivalent to $U$.  As seen above, the range of
$T_{\xi}^*T_{\xi}$ consists entirely of continuous functions.  But now, since
$T_{\xi}^*$ is an isometry, this range is a closed subspace of $L^2(G)$.  But
$T_{\xi}^*T_{\xi}$ is given by right convolution by $c_{\xi\xi}$.  We have thus
obtained:

\proclaim{8.7 Proposition} Let $U$ be a square-integrable irreducible
representation of $G$.  For any $\xi \in {\Cal B}_U$ normalized so that
$T_{\xi}^*$
is an isometry, right convolution by $c_{\xi\xi}$ is a projection of
$L^2(G)$ onto
a closed subspace consisting entirely of continuous functions, on which $L$ is
unitarily equivalent to $U$.
\endproclaim

Let $H_{\xi}$ denote the range of the isometry $T_{\xi}^*$.  (Note that
$c_{\xi\xi}
= T_{\xi}^*\xi$ so that $c_{\xi\xi} \in H_{\xi}$.)  For every $\varphi \in
H_{\xi}$ and every $x \in G$ we see above that
$$
\varphi(x) = \varphi * c_{\xi\xi}(x) = \int {\bar
c}_{\xi\xi}(x^{-1}y)\varphi(y)dy
= \<L_xc_{\xi\xi},\varphi\>.
$$
The second equality says exactly that $H_{\xi}$ is a ``reproducing-kernel
Hilbert
space'' on $G$, with reproducing kernel $c_{\xi\xi}$.  The third equality
says that
the map $x \mapsto L_xc_{\xi\xi}$ is a ``coherent state'' for $H_{\xi}$.  (See
\cite{Al} for a recent review of coherent states, with many 
interesting examples.)

Suppose now that $\zeta,\omega \in {\Cal B}_U$.  Then {}from $8.2$ but with
the roles
of the vectors interchanged, and {}from the definition of $\gamma(\zeta,\omega)$,
we obtain
$$
\align
\gamma(\omega,\zeta)\<\xi,\xi\> &= \<c_{\omega\xi},c_{\zeta\xi}\> 
= \int {\bar c}_{\omega\xi}(y)c_{\zeta\xi}(y)dy = \int
c_{\xi\omega}(y^{-1}){\bar c}_{\xi\zeta}(y^{-1})dy \\
&= \int {\bar c}_{\xi\zeta}(y)c_{\xi\omega}(y)\Delta(y^{-1})dy
= \<\Delta^{-1/2}c_{\xi\zeta},\Delta^{-1/2}c_{\xi\omega}\> 
= \<\Delta^{-1/2}T_{\xi}^*\zeta, \Delta^{-1/2}T_{\xi}\omega\>.
\endalign
$$
Thus we have
$$
\gamma(\omega,\zeta) =
\|\xi\|^{-2}\<\Delta^{-1/2}T_{\xi}^*\zeta,\Delta^{-1/2}T_{\xi}^*\omega\>.
$$
Notice that this is all well-defined, since as seen above,
$$
T_{\xi}^*\zeta = c_{\xi\zeta} \in L^2(G) \cap L^2(G,\Delta^{-1}dx)
$$
so that $\Delta^{-1/2}c_{\xi\zeta} \in L^2(G)$.  In conclusion, we obtain:

\proclaim{8.8 Theorem} Let $U$ be a square-integrable irreducible
representation of
$G$.  Let $\xi \in {\Cal B}_U$, normalized so that $T_{\xi}^*$ is an isometry.
Then for $\eta,\zeta \in {\Cal B}_U$ we have
$$
\int \alpha_x(\<\eta,\zeta\>_K)dx =
\|\xi\|^{-2}\<\Delta^{-1/2}T_{\xi}^*\zeta,\Delta^{-1/2}T_{\xi}^*\eta\> I_{H}.
$$
\endproclaim

We remark that if we choose $\eta$ and $\zeta$ so that
$$
\<\Delta^{-1/2}T_{\xi}^*\zeta,\Delta^{-1/2}T_{\xi}^*\eta\> = \|\xi\|^2,
$$
then we obtain
$$
\int \alpha_x(\<\eta,\zeta\>_K)dx = I_{H}.
$$
This is just another way of writing the familiar ``resolution of the identity''
{}from the theory of coherent states.

If $G$ is unimodular, we see that the right-hand-side of the equation of
Theorem
$8.8$ simplifies to
$$
\|\xi\|^{-2} \<\zeta,\eta\>I_{H},
$$
and now $\|\xi\|^2$ is the familiar formal dimension of $U$.  (See
\cite{D,Rf1}.)

If $G$ is not unimodular the right-hand-side of the equation of Theorem
$8.8$ is a
bit unattractive because the vectors and inner-product of the right-hand
side are
taken in $L^2(G)$, not $H$.  But the considerations just before the
statement of Theorem 8.8 show that if $\eta \in {\Cal B}_U$ then
$\Delta^{-1}T_\xi^*\eta \in L^2(G)$. Thus we can define an unbounded
operator, $K$, with dense domain ${\Cal B}_U$ by
$$
  K^{-1}\eta = \|\xi\|^{-2}T_{\xi}\Delta^{-1}T_{\xi}^*\eta.
$$
Then one can check that $K$ is a positive operator, and
the right-hand side of the equation of Theorem $8.8$ can be rewritten as
$$
  \<\zeta, K\eta\>I_H = \<K^{-1/2}\zeta,K^{-1/2}\eta\> I_{H}.
$$
This is the form given in theorem $3$ of \cite{DM}, or theorem $4.3$ of
\cite{Ca}, or theorem 2 of \cite{BT}.
We omit the details about domains and the fact that $K$ is independent of the
choice of $\xi$.  But one can check that, as expected {}from \cite{DM}, $K$ is
$\Delta^{-1}$-semi-invariant, reflecting the situation for the left regular
representation seen earlier.

We remark that in \cite{Mo} Moore has given orthogonality relations for factor
square-integrable representations.  But his orthogonality relations are not
for the
coefficient functions as defined here.  So it is not clear to me how his
results
relate to those given here.

We conclude this section by showing that the possible difficulty 
mentioned before Definition 1.2, namely that it may happen that 
$a \in {\Cal M}$  but $|a| \notin {\Cal M}$, actually occurs even in the 
present setting of square-integrable representations.

Let $G$ be the ``$ax+b$''  group. So $G$ 
is ${\Bbb R}\times {\Bbb R}^+$ with product given by
$$
    (p,s)(q,t) = (p+sq, st).
$$
The modular function of $G$ is $\Delta(p,s) = s^{-1}$. It is 
well-known that $G$ has two inequivalent square-integrable 
irreducible representations. We consider one of them. It has many 
models. For our purposes the most elementary approach to what 
we need seems to be given by theorem 2 of \cite{BT}, so we use the 
model used there. The Hilbert space is $H = L^2({\Bbb R}^+, dt/t)$, 
and the representation is given by
$$
  (\pi_{(p,s)}\xi)(t) = e(pt)\xi(st),
$$
where by definition $e(t) = \exp(2\pi it)$. Let $K$ denote the 
unbounded operator on $H$ defined by
$$
  (K\xi)(t) = t\xi(t).
$$
One can check that $K$ is $\Delta^{-1}$-invariant.

Applying theorem 2 of \cite{BT} and Theorem 8.8 above to this 
particular situation, we find that, up to multiplication by a positive 
scalar, 
$$
  \int_G \a_x(T) \ dx = Tr(K^{-1}T)I_H
$$
for any positive compact operator $T$, where $\a$ is the action of 
conjugation by $\pi$. As in our earlier discussion of the irreducible 
case, this gives essentially an ordinary weight. We can now study 
this weight independently of the fact that it comes {}from a group 
representation.

To simplify our analysis, we make the change of variables 
$r = e^{-t}$.  Then our Hilbert space becomes $L^2({\Bbb R})$ for 
Lebesgue measure, and $D = K^{-1}$ is the operator of pointwise 
multiplication by $t \mapsto e^t$. We denote our weight by $\psi$. 
It is now given by 
$$
  \psi(T) = Tr(DT)   .
$$
Since $D$ is unbounded, we must make precise what this means. 
We can do this conveniently in terms of spectral projections of $D$. 
For our purposes the following works well. For each integer $n \geq 
1$ let $E_n$ denote the orthogonal projection onto the subspace of 
$L^2({\Bbb R})$ consisting of the functions supported in 
$[-n,\ -n+1]\cup[n-1,\ n]$. Thus the $E_n$'s are mutually 
orthogonal and sum to $I$. For each $n$ the operator $DE_n$ is a 
bounded positive operator. For any positive bounded operator $T$ 
we take $\psi(T)$ to mean
$$
  \psi(T) = \sum Tr(DE_nT)  .
$$
In particular, $\psi$ is lower semi-continuous. We 
let ${\Cal M}^+ = \{T \in B(H)^+: \psi(T) < \infty\}$,  and we 
let $\Cal M$ be the linear span of ${\Cal M}^+$. For our discussion of 
square-integrable representations we are most interested in the 
restriction of $\psi$ to the algebra $K(H)$ of compact operators.

\proclaim{8.9 Theorem} With notation as just above, there are $S, 
T \in K(H)$ such that $S, T \in {\Cal M}^+$ but $|S-T| \notin {\Cal M}^+$.
\endproclaim

\demo{Proof} For each integer $n \geq 1$ 
set $\xi_n = \chi_{[n-1, \ n]}$, and $\eta_n = \chi_{[-n,\  -n+1]}$, 
where $\chi$ denotes ``characteristic function''. Thus $\xi_n$ and 
$\eta_n$ are unit vectors in the range of $E_n$. The following steps 
are motivated by the example following theorem 2.4 of \cite{Pe1}. 
Choose a sequence $\{a_n\}$ of numbers with $0 \leq a_n \leq 1$ 
such that $\sum a_n < \infty$ but $\sum a_n^{1/2} = \infty$. For 
instance, $a_n = n^{-2}$. Let $P_n$ and $Q_n$ denote the rank--1 
projections which are  $0$ on the orthogonal complement of 
$\{\xi_n, \eta_n\}$, whereas with respect to the 
basis $\{\xi_n, \eta_n\}$ they have matrices
$$
\left( \matrix
a_n & (a_n - a_n^2)^{1/2} \\
(a_n - a_n^2)^{1/2} & 1 - a_n
\endmatrix \right)
$$ 
and 
$$
\left( \matrix
0 & 0 \\
0 & 1
\endmatrix \right)   ,
$$
respectively. Then $(P_n \ -\ Q_n)^2$ has matrix
$$\left( \matrix
a_n & 0 \\
0 & a_n
\endmatrix \right)  ,
$$
so that $|P_n\  -\ Q_n|$ has matrix
$$
\left( \matrix
a_n^{1/2} & 0 \\
0 & a_n^{1/2}
\endmatrix \right)   .
$$

For any integer $k$ with $|k| \geq 1$ set $d_k = \<D\xi_k, \xi_k\>$ if 
$k\geq 1$ and $d_k = \<D\eta_{-k}, \eta_{-k}\>$ if $k \leq -1$. Thus, 
disregarding $k = 0$, we see that $\{d_n\}$ goes to $0$ rapidly as $k 
\to -\infty$, and to $+\infty$ rapidly as $k \to +\infty$. We can use 
the basis $\{\xi_n, \eta_n\}$ to evaluate $\psi$ on $P_n$ and $Q_n$, 
and a quick calculation shows that
$$
  \psi(P_n) = a_nd_n\ +\ (1-a_n)d_{-n},
$$
$$
   \psi(Q_n) = d_{-n},
$$
$$
   \psi(|P_n - Q_n|) = a_n^{1/2}(d_n\ +\ d_{-n})   .
$$

Set $S = \sum d_n^{-1}P_n$ and $T = \sum d_n^{-1}Q_n$, where the 
sums are for $n \geq 1$. Since $\{d_n\}$ grows rapidly for positive 
$n$, the sums converge in norm, and $S, T \in K(H)^+$. {}From the 
above calculations and the properties which we required of 
$\{a_n\}$ we see that
$$
   \psi(S) = \sum a_n \ + \ \sum d_n^{-1}(1-a_n)d_{-n} \  < \infty  ,
$$
$$
  \psi(T) = \sum d_n^{-1}d_{-n} \ < \ \infty.
$$
Thus $S, T \in {\Cal M}^+$. However
$$
  \psi(|S-T|) = \sum a_n^{1/2} \ + \ \sum a_n^{1/2}d_n^{-1}d_{-n} \ 
= \ \infty   .
$$
Thus$|S-T| \notin {\Cal M}^+ .$
\qed
\enddemo

\enddocument